\documentclass[a4paper,12pt]{article}

\newenvironment{proof}{\hbox{}\vspace{-0.5cm} {\bf Proof:}}{\hfill $\Box$ \\}

\usepackage{lmodern}

\newtheorem{theorem}{Theorem}
\newtheorem{lemma}{Lemma}
\newtheorem{proposition}{Proposition}

\newtheorem{algorithm}{Algorithm}

\newtheorem{assumption}{Assumption}
\newtheorem{remark}{Remark}

\newtheorem{example}{Example}

\usepackage{amsmath}
\usepackage{amssymb}
\usepackage{amsfonts}
\usepackage{mathrsfs} 
\usepackage{graphicx}
\usepackage{color}
\usepackage{ulem}
\usepackage{bm}
\usepackage{booktabs}
\usepackage{algorithm}
\usepackage{algorithmic}
\usepackage{hyperref}
\newcommand{\R}{\mathbb{R}} % real numbers
\renewcommand{\S}{\mathbb{S}} % symmetric matrices
\newcommand{\trace}{\mathrm{tr}}

\newcommand{\ba}{\mathbf{a}}
\newcommand{\bb}{\mathbf{b}}

\newcommand{\bp}{\mathbf{p}}
\newcommand{\bq}{\mathbf{q}}
\newcommand{\bu}{\mathbf{u}}
\newcommand{\bv}{\mathbf{v}}
\newcommand{\bx}{\mathbf{x}}
\newcommand{\bw}{\mathbf{w}}
\newcommand{\by}{\mathbf{y}}
\newcommand{\bz}{\mathbf{z}}
\newcommand{\bmu}{\boldsymbol{\mu}}
\DeclareMathOperator{\Tr}{tr}

\DeclareMathOperator{\Int}{int}

\DeclareMathOperator{\conv}{conv}

\title{\bf Maximal entropy in the moment body}

\begin{document}

\author{Didier Henrion$^{1,2}$}

\footnotetext[1]{CNRS; LAAS; Universit\'e de Toulouse, 7 avenue du colonel Roche, F-31400 Toulouse, France.}
\footnotetext[2]{Faculty of Electrical Engineering, Czech Technical University in Prague,
Technick\'a 2, CZ-16626 Prague, Czechia.}

\date{Draft of \today}

\maketitle
\addtocounter{footnote}{2}

\begin{abstract}
	A moment body is a linear projection of the spectraplex, the convex set of trace-one positive
	semidefinite matrices. Determining whether a given point lies within a given moment body is a problem
	with numerous applications in quantum state estimation and polynomial optimization. This moment body
	membership oracle can be addressed with semidefinite programming, for which several off-the-shelf
	interior-point solvers are available. In this paper, inspired by techniques from quantum information
	theory, we argue analytically and geometrically that a much more efficient approach consists of
	minimizing globally a smooth strictly convex log-partition function, dual to a maximum entropy
	problem. We analyze the curvature properties of this function, showing that conditioning is governed
	by the distance of the point to the boundary of the moment body, and we describe a neat geometric
	preconditioning algorithm that exploits this analysis.  Basic numerical experiments, comparing against interior-point and first-order
	semidefinite solvers, reveal a cubic dependence on the matrix size, similar   to a few eigenstructure
	computations.  They also illustrate the two regimes of the oracle: dense projections are handled
	efficiently up to sizes of several hundred, while sparse instances such as matrix completion scale to
	matrices of size several thousand in minutes on a standard laptop. In both cases the main bottleneck
	in this approach to large-scale semidefinite programming is moved almost entirely to efficient gradient storage and
	manipulation.
\end{abstract}

\section{Introduction}

Semidefinite programming is a versatile framework for convex optimization. It consists of optimizing
(typically linear functions) over spectrahedra (linear sections of the semidefinite cone, described by
linear matrix inequalities) or spectrahedral shadows (linear projections of spectrahedra). These sets
capture a large class of convex semialgebraic sets \cite{bn01}. Polynomial optimization relies heavily
on semidefinite optimization, and the moment-SOS hierarchy constructs a nested family of spectrahedral
shadows of increasing size that provide increasingly tight approximations of convex hulls of
semialgebraic sets, see e.g. \cite{hkl20,n23,t24} and references therein.

Semidefinite optimization problems can be solved with interior-point algorithms \cite{nn94,bn01}.
However, as second-order methods, these algorithms do not scale well at the age of data science. Most
of the computational burden is concentrated on computing and storing the Hessian matrix of
second-order derivatives of a logarithmic barrier function. First-order algorithms scale better, since
they use only gradient information, but they are also more sensitive to problem scaling and
conditioning. Conditioning of semidefinite optimization problems is understood theoretically
\cite{renegar95}, but evaluating the conditioning of a given problem is as expensive as solving the
original problem. From that point of view, the versatility and generality of semidefinite programming
can also be seen as a weakness: currently, there is no simple recipe that can be systematically used to
cure all numerical issues, see e.g. \cite{pataki19} for a survey of recent attempts. There are at
least three geometric pathologies that can occur in semidefinite programming: (i) a linear image of an
unbounded spectrahedron need not be closed; (ii) a spectrahedron or its shadow can lack interior
points; (iii) the linear map defining a spectrahedral shadow can be ill-conditioned (i.e. with singular
values largely differing in magnitude). In this paper, we propose to focus on pathology (iii), namely
ill-conditioning of the linear map, and our strategy is as follows. First, we restrict our attention to
semidefinite feasibility problems whose spectrahedral shadows are full-dimensional and bounded. This
eliminates the pathologies (i) and (ii). Second, we focus on analytic, quantitative aspects of a
standard first-order optimization algorithm, in which issue (iii) appears explicitly through curvature
parameters. This allows us to design a simple and cheap pre-conditioning algorithm.

Our focus is on the moment body membership oracle problem: finding a point in a linear projection of
the spectraplex, defined as the compact convex set of trace-one positive semidefinite matrices, a
non-polyhedral generalization of the simplex. Determining whether a given point lies within a given
moment body is a problem with numerous applications in polynomial optimization or quantum information
theory. This includes for example the problem of decomposing a given multivariate polynomial as a sum
of squares (SOS) of other polynomials, see e.g. \cite[Section 2.4]{n23} and references therein. In
order to address this problem with a first-order algorithm, we use an approach inspired from quantum
information theory \cite{jv23,h24}, namely the global minimization of a smooth and strictly convex
log-partition function dual to a maximum entropy problem. Quantum state estimation aims to recover a
density matrix (i.e. an element of the spectraplex, a trace-one positive semidefinite matrix)
consistent with observed measurement statistics (i.e. the linear projection of the spectraplex)
\cite{bz17} --- and this is exactly our moment body membership oracle problem. A particularly effective
method for solving this problem consists of selecting, among all compatible density matrices, the one
maximizing entropy. The dual of this problem leads to the minimization of a convex, smooth function
called the log-partition function. We analyze its curvature properties and, based on geometric
quantities appearing during the analysis, we describe a neat and simple geometric pre-conditioning
algorithm; the same dual viewpoint also yields exact, rational certificates separating strict
feasibility, boundary membership and infeasibility.

Our approach belongs to a recent line of work that solves semidefinite programs through the dual of an
entropic regularization. Reference \cite{l23} regularizes a general SDP with the von Neumann entropy (or
an unnormalized variant), solves the resulting dual, and estimates the dual gradients with randomized
trace estimators, thereby avoiding the cubic cost of a full spectral factorization; the method is
specialized to SDPs with a diagonal constraint and to spectral projector problems, with applications to
combinatorial optimization and spectral embedding. The follow-up reference \cite{cl25b} develops this into a
non-Euclidean dual gradient ascent, in which the dual is maximized with respect to a problem-adapted
norm; they prove that the dual gradient norm converges to zero at a rate independent of the ambient
dimension, recover primal-feasible points by rounding in special cases, and again rely on randomized
trace estimators for scalability, with the MaxCut, optimal transport and permutation
synchronization relaxations as guiding examples. The present work shares this entropic-dual standpoint
but differs in scope and in several theoretical respects. We do not regularize a general objective: we
treat the membership problem on the spectraplex, where the trace-one normalization is exact rather than
introduced as a regularizer. The dual is then precisely the log-partition function $\log\trace\exp$,
which is globally smooth and strictly convex, in contrast with the partition function $\trace\exp$
arising from the unnormalized regularization, whose gradient has a much larger dynamic range. This
regularity is what lets us carry out an explicit Euclidean curvature analysis: a simple and cheap
whitening of the data makes the function $\tfrac12$-smooth, and on the relevant sublevel set its
strong-convexity modulus, and hence the condition number  is controlled by the geometric feasibility margin. Whereas \cite{cl25b} absorbs the
problem geometry into an adapted norm and reports a dimension-independent rate, our analysis makes the
dependence explicit and locates the only source of ill-conditioning at the boundary of the moment body. We use exact gradients and off-the-shelf L-BFGS to reach accuracy $10^{-8}$, and we
give a complete treatment of the decision problem through the exact certificates mentioned above. The
two viewpoints are complementary: the randomized trace estimators of \cite{l23,cl25b} are precisely the
tool that could carry our preconditioned dual past its current per-iteration spectral factorization,
the main scalability limit in the dense regime.

Basic numerical experiments with a rudimentary Matlab prototype illustrate the approach in two
complementary regimes. It is one to two orders of magnitude faster than the interior-point solver MOSEK \cite{mosek}
and competitive with SDPNAL+ \cite{sdpnal15,sdpnal20}, a state-of-the-art first-order solver for
large-scale semidefinite programming. On dense (i.e. non-sparse) projections, where the one-off
pre-conditioning dominates the cost, problems of size up to several hundred are solved at expected
accuracy $10^{-8}$ in a fraction of a second; here the limiting factor is the storage and manipulation
of the dense linear map and iterate. On sparse instances such as positive semidefinite matrix
completion, where the data is orthonormal by construction and pre-conditioning is essentially free, the
method scales much further: feasible completions of matrices of size up to $n=5000$, with millions of
revealed entries, are solved to the same accuracy in minutes on a standard laptop, and a problem of
size $n=1000$ in about a second. Practically speaking, this implies that, for this problem class, the
bottleneck of this class of large-scale semidefinite solvers is pushed almost exclusively to the efficient storage
and manipulation of gradient information, consistently with the randomized and sketching techniques
of \cite{l23,cl25,sketchy,hallar}.

\subsection{Outline}

The paper is organized as follows. Section~2 defines the moment body and presents a few examples to
illustrate its geometry in low dimension. In Section~3 we show how testing membership in a moment body
of size $m$ defined by a spectraplex of size $n$-by-$n$ can be formulated as the unconstrained
minimization in $\R^m$ of a smooth, strictly convex log-partition function, and we prove that this dual
problem is equivalent (via strong duality) to a primal maximum entropy formulation. Section~4 is
devoted to a first geometric analysis of the dual objective. We derive explicit upper and lower bounds
on its Hessian in terms of the centered Gram matrix of the linear map defining the moment body, showing
that the dual is globally $\lambda$-smooth and $\alpha$-strongly convex on sublevel sets, with $\lambda$
and $\alpha$ depending on the spectrum of that Gram matrix and on the distance of the queried point to
the boundary of the moment body. In Section~5 we present a simple preconditioning algorithm: by
centering and orthonormalizing the linear map, one forces the dual to become $1/2$-smooth, while
the strong-convexity modulus on the relevant sublevel set becomes $1/n^{1+\sqrt2/\delta}$, where
$\delta$ is the distance of the queried point to the boundary. Therefore the conditioning is polynomial
in $n$ for a fixed margin and degrades only as the point approaches the boundary. In Section~6 we
exploit these curvature estimates to bound the value and the norm of the unique minimizer in terms of
the input data. Section~7 gives a detailed iteration complexity analysis of a gradient algorithm
applied to the preconditioned dual. Section~8 discusses how the same dual framework decides the three
cases of the membership problem, providing computable certificates of strict feasibility,
boundary membership (weak feasibility) and infeasibility. Section~9 briefly explains how our analysis
extends to block-separable (direct-sum) moment-body problems, in which the primal density matrix splits
into several independent blocks. Finally, Section~10 presents numerical experiments on random instances
in two complementary regimes: we compare our Matlab prototype against off-the-shelf semidefinite
solvers, solving dense problems up to size of several hundred and sparse instances such as matrix
completion up to size $n=5000$ on a standard laptop.

\subsection{Notations}

$\S^n$ is the space of real valued symmetric matrices of size $n$, $$\S^n_+ := \{X \in \S^n, X \succeq 0\}$$ is the convex closed cone of positive semidefinite elements of $\S^n$, called the {\it semidefinite cone}. Its interior $$\mathrm{int}\:\S^n_+ := \{X \in \S^n, X \succ 0\}$$ is the convex open cone of positive definite elements of $\S^n$, and $$\S^n_1:=\{X \in \S^n, X \succeq 0, \trace X = 1\}$$ is called the {\it spectraplex}, a generalization to non-diagonal matrices of the polyhedral simplex. It is a spectrahedron, an affine slice of the semidefinite cone, see {e.g. \cite{BPT13} or} \cite[Section 7.3]{n23}. 

Given a {symmetric} matrix $X$, $\log X$ denotes its logarithm, $\exp X$ its exponential, and $$\exp_1 X := \frac{\exp X}{\trace \exp X}$$ is the normalized or trace-one exponential.
 
\section{The moment body}

Let $A_i \in \S^n$, $i=1,\ldots,m$ be given matrices.
Define the linear map $\mathcal A : \S^n \to \R^m, X \mapsto \trace(A_i X)_{i=1,\ldots,m}$ and its adjoint $\mathcal A^T : \R^m \to \S^n, \by \mapsto A(\by):=\sum_{i=1}^m y_i A_i$.
The {\it moment body} of $\mathcal A$ is the set
\[
\boxed{
{\mathscr M}:=  \left\{ \mathcal A(X) : X \succeq 0, \trace X = 1\right\} \subset \R^m}
\]
or equivalently $\mathscr M:=\mathcal A(\S^n_1)$. In words, a moment body is the linear image of a spectraplex.
As a linear projection of a convex and compact set, set $\mathscr M$ is also convex and compact. The terminology moment body is motivated as follows. 
Let $\mathscr X$ be a {set in a topological vector space}, and let $\phi : \mathscr X \to \mathscr U$ be a surjective map, {where $\mathscr U:=\{\bu \in \R^n : \bu^T \bu = 1\}$ is the unit sphere.}
For example, $\phi(\bx)$ can be the result of a measurement for $\bx \in \mathscr X$ in some given set of Euclidean space,  {satisfying the normalization identity
$\phi^T(\bx)\phi(\bx)=1$ for all $\bx \in \mathscr X$.}
We can identify each matrix $A_i$ with the Gram matrix of a function $a_i : \mathscr X \to \R, {\bx} \mapsto \phi^T(\bx) A_i \phi(\bx)$ and then write $\trace(A_iX) = \int_{\mathscr X} a_i(\bx) d\mu(\bx)$ where $X = \int_{\mathscr X} \phi(\bx) \phi(\bx)^T d\mu(\bx)$ is the moment matrix of $\mu$, an element of $\mathscr P(\mathscr X)$, the set of probability measures on
$\mathscr X$. Equivalently, if we define $\nu$ as the image measure of $\mu$ through $\phi$, $X = \int_{\mathscr U} \bu\bu^T d\nu(\bu)$ is the {autocorrelation} matrix of $\nu \in \mathscr P(\mathscr U)$.
Both measures satisfy $\int_{\mathscr X} d\mu(\bx) = \int_{\mathscr U} d\nu(\bu) = \int_{\mathscr U} \bu^T \bu \,d\nu(\bu) = \trace X =1$.
The moment body is therefore the set of all moments of such probability measures, i.e.
\[
\mathscr M = \left\{ \int_{\mathscr X} \ba(\bx) d\mu(\bx) \::\: \mu \in \mathscr P(\mathscr X)\right\}.
\]
{If $X$ is complex Hermitian, $\S^n_1$} is also called the set of mixed quantum states in quantum information theory \cite{bz17}. Its elements are known as density operators or density matrices. Its extreme points are rank-one matrices generated by {vectors of the complex unit sphere.} Alternatively, we can also interpret the moment body as a {joint} numerical range -- see e.g. \cite{bd71,psw21,nv24} and references therein --
defined as the convex hull of the image of the {complex} unit sphere through the linear map $\mathcal A$. {Note however that in this paper we restrict our attention to real matrices, and the moment body can be interpreted as} the convex hull of the image of the unit sphere through the linear map, i.e.
\[
\mathscr M = \conv \mathcal A(\mathscr U) = \conv\{[\bu^T A_i \bu]_{i=1,\ldots,m}, \bu \in \mathscr U\}.
\]
Finally, as a linear projection of a spectrahedron, the moment body is a {\it spectrahedral shadow}, see e.g. \cite{BPT13} or \cite[Section 7.3]{n23}. 

{	
	Note however that not all bounded spectrahedra or spectrahedral shadows can be modeled as moment bodies. Given a vector $\by \in \R^m$, the support function of $\mathscr M$ is
\begin{equation}\label{eq:support}
h(\by) := \max_{\bz \in \mathscr M} \by^T \bz = \max_{X \in \S^n_1} \Tr(A(\by)X) = \min_{\lambda \in \R} \{\lambda : \lambda I_n - A(\by) \succeq 0\}
\end{equation}
where the last identity follows from strong duality and strict feasibility of $X = \tfrac{1}{n} I_n$. Therefore the polar body
\[
\mathscr M^o:=\{\by \in \R^m : h(\by) \leq 1\}=\{\by \in \R^m : I_n - A(\by) \succeq 0\}
\]
must be a spectrahedron.
For example, the elliptope
\[
\left\{\bz \in \R^3 : \left(\begin{array}{ccc}1 & z_1 & z_2 \\ z_1 & 1 & z_3 \\ z_2 & z_3 & 1\end{array}\right) \succeq 0\right\}
\] is not a moment body because its polar, the convex hull of the Roman Steiner surface, is not a spectrahedron, see e.g. \cite[Ex. 5.4.4]{BPT13}. This polar is a moment body by construction.}

Throughout the paper, we make the following natural assumption on the linear map.

\begin{assumption}\label{injective}
{Matrices $I_n, A_1, \ldots, A_m$ are linearly independent in $\S^n$.}
\end{assumption}

{Note that Assumption \ref{injective} implies that linear map $\mathcal A$ is surjective. It also implies that $\mathcal A(\tfrac{1}{n} I_n)$ is an interior point of $\mathscr M$.}

\begin{figure}[h]
\begin{center}
	\includegraphics[width=0.7\textwidth]{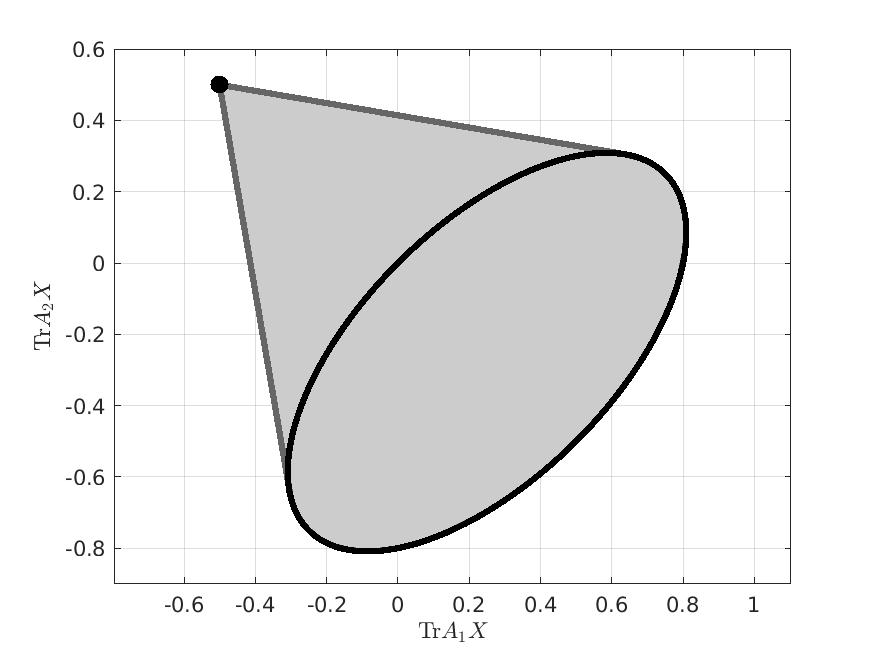}
	\caption{The moment body (light gray) of Example \ref{ex:momentbody2d} is the convex hull of an ellipse (black, bottom right) and a point (black, top left).\label{fig:momentbody2d}}
\end{center}
\end{figure}
\begin{example}\label{ex:momentbody2d}
Let $n=3$, $m=2$ and
\[
A_1 = 
\frac{1}{2}\left(\begin{array}{rrr}
	1 & 1 & 0 \\ 1 & 0 & 0 \\ 0 & 0 & -1
\end{array}\right), \quad
A_2 = 
\frac{1}{2}\left(\begin{array}{rrr}
-1 & 1 & 0 \\ 1 & 0 & 0 \\ 0 & 0 & 1
\end{array}\right).
\]
As explained e.g. in \cite{h10}, the moment body of $\mathcal A$ is the convex hull of the algebraic curve dual to the curve $\{ \by \in \R^2 : p(\by)=\det(I_3+A_1y_1+A_2y_2)=0\}$, i.e. the envelope of all tangent lines. The determinant factors into $p(\by)=\frac{1}{8}(4+2y_1-2y_2-y^2_1-2y_1y_2-y^2_2)(2-y_1+y_2)$, so the dual curve is the union of the ellipse $\{\bx \in \R^2 : 5x^2_1-6x_1x_2+5x^2_2-4x_1+4x_2=0\}$ and the point $(-\frac{1}{2},\frac{1}{2})$. Equivalently, in parametric form, the moment body of $\mathcal A$ is the convex hull of the ellipse
\[
\left\{\left(\trace\frac{1}{2}\left(\begin{array}{ccc}1 & 1 & 0 \\ 1 & 0 & 0 \\ 0 & 0 & -1 \end{array}\right)
\left(\begin{array}{c}\cos\theta \\ \sin\theta \\ 0\end{array}\right)
\left(\begin{array}{c}\cos\theta \\ \sin\theta \\ 0\end{array}\right)^T,
\trace\frac{1}{2}\left(\begin{array}{ccc}-1 & 1 & 0 \\ 1 & 0 & 0 \\ 0 & 0 & 1 \end{array}\right)\left(\begin{array}{c}\cos\theta \\ \sin\theta \\ 0 \end{array}\right)
\left(\begin{array}{c}\cos\theta \\ \sin\theta \\ 0 \end{array}\right)^T\right),  \theta \in [0,2\pi]\right\}
\]
and the point
\[
\left\{\left(\trace\frac{1}{2}
\left(\begin{array}{ccc}1 & 1 & 0 \\ 1 & 0 & 0 \\ 0 & 0 & -1 \end{array}\right)
\left(\begin{array}{c}0 \\ 0 \\ 1 \end{array}\right)
\left(\begin{array}{c}0 \\ 0  \\ 1 \end{array}\right)^T,
\trace\frac{1}{2}\left(\begin{array}{ccc}-1 & 1 & 0 \\ 1 & 0 & 0 \\ 0 & 0 & 1 \end{array}\right)
\left(\begin{array}{c}0 \\ 0 \\ 1 \end{array}\right)
\left(\begin{array}{c}0 \\ 0  \\ 1 \end{array}\right)
\right) \right\}.
\] 
See Figure \ref{fig:momentbody2d}.
\end{example}

%\begin{example}
%	Let $n=m=3$ and
%\[
%A_1 = \begin{pmatrix}0 & 1 & 0\\1 & 0 & 0\\0 & 0 & 0\end{pmatrix}, 
%A_2 = \begin{pmatrix}0 & 0 & 1\\0 & 0 & 0\\1 & 0 & 0\end{pmatrix}, 
%A_3 = \begin{pmatrix}0 & 0 & 0\\0 & 0 & 1\\0 & 1 & 0\end{pmatrix}.
%\]
%
%steiner's roman surface - dual to cayley
%\end{example}

\begin{figure}[h]
	\begin{center}
		\includegraphics[width=0.9\textwidth]{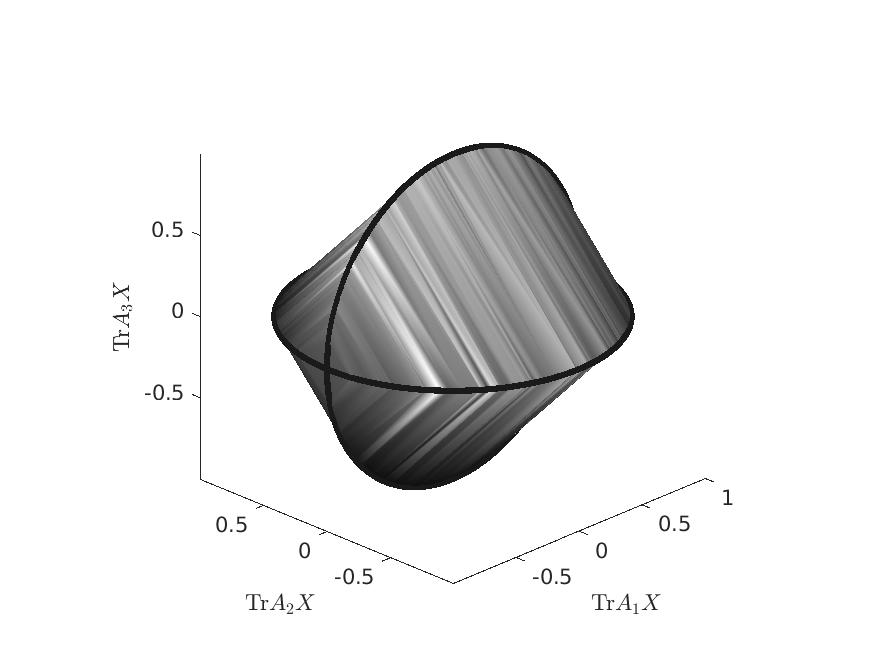}
		\caption{The moment body (light gray) of Example \ref{ex:momentbody3d} is the convex hull of two orthogonal circles (black).\label{fig:momentbody3d}}
	\end{center}
\end{figure}
\begin{example}\label{ex:momentbody3d}
Consider the two unit circles in orthogonal planes in~$\R^3$:
\[
C_{x_1x_2} = \{(\cos\theta,\;\sin\theta,\;0)\colon\theta\in[0,2\pi]\}, 
\quad
C_{x_1x_3} = \{(\cos\phi,\;0,\;\sin\phi)\colon\phi\in[0,2\pi]\}.
\]
Their convex hull can be modeled as a moment body as follows. Define the matrices
\[
J = 
\begin{pmatrix}
1 & 0 \\
0 & -1
\end{pmatrix},
\:
K = 
\begin{pmatrix}
0 & 1 \\
1 & 0
\end{pmatrix},
\:
A_1 = \begin{pmatrix} J & 0 \\[3pt] 0 & J \end{pmatrix},
\:
A_2 = \begin{pmatrix} K & 0 \\[3pt] 0 & 0_{2} \end{pmatrix},
\:
A_3 = \begin{pmatrix} 0_{2} & 0 \\[3pt] 0 & K \end{pmatrix},
\]
i. e. $A_i \in \S^4$, $i=1,2,3$.
The corresponding moment body is the convex hull of the union of the circles $C_{x_1x_2}$ and $C_{x_1x_3}$. Indeed, $\{\mathcal A(X_1\oplus0_{2})  : X_1 \in \S^2_1\} = \{(\trace X_1 J), \trace (X_1 K), 0) : X_1 \in \S^2_1\} = \conv C_{x_1x_2}$, $\{\mathcal A(0_2\oplus X_2)  : X_2 \in \S^2_1\} =\{(\trace X_2 J), 0, \trace (X_2 K)) : X_2 \in \S^2_1\} = \conv C_{x_1{x_3}}$, and $\mathcal A(\S^4_1)$ consists of all convex combinations of these two sets, see Figure \ref{fig:momentbody3d}.
\end{example}

\section{The moment body membership oracle}

Given the linear map $\mathcal A$, the moment body membership oracle consists of determining whether a given vector $\bb \in \R^m$ belongs to $\mathscr M$. 

Let
\[
\boxed{f(\by):=\log\trace\exp A(\by) - \bb^T \by}
\]
{where $\log\trace\exp$ is} the cumulant generating function or log-partition function, {a standard terminology in information geometry and quantum statistical mechanics.}

\begin{lemma}\label{gh}
Function $f$ is smooth and convex on $\R^m$. Its gradient is $$\nabla f(\by) = \left[\trace(A_i X(\by))-b_i\right]_{i=1,\ldots,m} = \mathcal A(X(\by))-\bb$$ and its Hessian is	
\[{
\nabla^2 f(\by) =\left[ \int_{0}^{1}
		\trace \bigl(A_i\, X(\by)^s\,A_j\,X(\by)^{1-s}\bigr)\,ds
	{ - } \trace (A_i\,X(\by) )\,\trace (A_j\,X(\by))\right]_{i,j=1,\ldots,m}.
}\]
where $$X(\by):=\exp_1 A(\by) \in \S^n_1$$ is a so-called density matrix.  
\end{lemma}

	\begin{proof}
		By standard matrix‐calculus, \(X\mapsto \exp X\) is smooth on \(\S^n\), so
		\(\log\Tr\exp X\) is smooth as a composition. 
Let \(t(\by)=\trace\exp A(\by)\) so that \( f(\by)=\log t(\by)\;-\;\by^Tb.\)

\medskip\noindent{\bf First derivatives. }
{
By the Duhamel formula for the derivative of a matrix exponential \cite{Wilcox1967}, we have
\begin{equation}\label{duhamel}
\frac{\partial \exp A(\by)}{\partial y_i}
=\int_{0}^{1}
\exp\bigl((1-s)A(\by)\bigr)\,A_i\,\exp\bigl(sA(\by)\bigr)\,ds.
\end{equation}
Taking the trace gives
\begin{equation}\label{eq1}
\frac{\partial t(\by)}{\partial y_i}
=\int_{0}^{1}
\trace\bigl(\exp((1-s)A(\by))\,A_i\,\exp(sA(\by))\bigr)\,ds
=\trace\bigl(A_i\,\exp A(\by) \bigr)
\end{equation}
by the cyclic property of the trace. Therefore
\[
\frac{\partial \log t(\by)}{\partial y_i}
=\frac{1}{t(\by)}\,\frac{\partial t(\by)}{\partial y_i}
=\frac{\trace(A_i\,\exp A(\by))}{\trace(\exp A(\by))}
=\trace(A_i\,X(\by))
\]
and finally
\[
\frac{\partial f(\by)}{\partial y_i} 
=\trace(A_i\,X(\by))-b_i.
\]
}

\medskip\noindent{\bf Second derivatives. }
{
Let us differentiate the gradient
\begin{equation}\label{eq2}
\frac{\partial^2 f(\by)}{\partial y_i \partial y_j} 
=\frac{\partial\,\trace(A_i\,X(\by))}{\partial y_j}\
=\trace \Bigl(A_i\,\frac{\partial X(\by)}{\partial y_j}\Bigr).
\end{equation}
First develop
\[
\frac{\partial X(\by)}{\partial y_j}
=\frac{1}{t(\by)}\,\frac{\partial \exp A(\by)}{\partial y_j}
\;-\;\frac{\exp A(\by)}{t(\by)^2}\,\frac{\partial t(\by)}{\partial y_j}
\]
and use relation \eqref{eq1}:
\[
\frac{\partial t(\by)}{\partial y_j}
=\trace(A_j \exp A(\by))
=t(\by)\,\trace(A_j X(\by))
\]
to obtain
\[
\frac{\partial X(\by)}{\partial y_j}
=\frac{1}{t(\by)}\,\frac{\partial \exp A(\by)}{\partial y_j}
\;-\;X(\by)\,\trace(A_j X(\by)).
\]
We use again Duhamel's formula \eqref{duhamel} to obtain
\[
\frac{\partial X(\by)}{\partial y_j}
=\int_{0}^{1}\frac{1}{t(\by)}\exp\bigl(s\,A(\by)\bigr)\,A_j\,\exp\bigl((1-s)\,A(\by)\bigr)\,ds
\;-\;X(\by)\,\trace\bigl(A_j\,X(\by)\bigr).
\]
Substituting this expression into relation \eqref{eq2} we get
\[
\frac{\partial^2 f(\by)}{\partial y_i \partial y_j} 
= \int_{0}^{1}\frac{1}{t(\by)}
\trace \Bigl(A_i\,\exp(sA(\by))\,A_j\,\exp((1-s)A(\by))\Bigr)\,ds
\;-\;\trace(A_i\,X(\by))\,\trace(A_j\,X(\by)).
\]
Since $X(\by)$ is symmetric and $s \in [0,1]$ it holds
\[
X(\by)^s=\frac{\exp(sA(\by))}{t(\by)^s},
\quad
X(\by)^{1-s}=\frac{\exp((1-s)A(\by))}{t(\by)^{1-s}}
\]
and we have
\[
\exp(sA(\by))\,A_j\,\exp((1-s)A(\by))
=  {t(\by)}\,X(\by)^s\,A_j\,X(\by)^{1-s}.
\]
Substituting this expression under the integral, we finally obtain
\begin{equation}\label{eq3}
\frac{\partial^2 f(\by)}{\partial y_i \partial y_j} 
	=\int_{0}^{1}
	\trace \bigl(A_i\,X(\by)^s\,A_j\,X(\by)^{1-s}\bigr)\,ds
	\;-\;
	\trace (A_i\,X(\by) )\,\trace (A_j\,X(\by))
\end{equation}
which is the expected expression.
Note that these expressions are classical in quantum information theory, see e.g. \cite [Lem. VI]{w63}, the Bogoliubov-Kubo-Mori (BKM) inner product in \cite[Sect. 7.3]{an00} or \cite[Prop. 6.1]{w14}.
}

\medskip\noindent{\bf Convexity. }
Let us show that for any direction \(\bv\in\R^m\) and any vector $\by \in \R^m$ it holds
\[
\bv^T\nabla^2 f(\by)\,\bv \;\ge0.
\]
From relation \eqref{eq3} it holds
\begin{equation}\label{eq4}
	\bv^T\nabla^2 f(\by)\,\bv
	=\int_0^1
	\Tr\bigl(A(\bv)\,X(\by)^s\,A(\bv)\,X(\by)^{1-s}\bigr)\,ds
	\;-\;
	\bigl[\Tr\bigl(A(\bv)\,X(\by)\bigr)\bigr]^2.
\end{equation}
For each fixed \(s\in[0,1]\), define the bilinear form
\[
\langle Y,Z\rangle_s \;:=\;\Tr\bigl(Y\, X^s(\by)\,Z\,X^{1-s}(\by)\bigr).
\]
This form satisfies the properties of an inner product {on $\S^n$} because $X(\by)$ is positive definite. 
Then by the usual Cauchy-Schwarz inequality for this inner product,
\[
\langle A(\bv),A(\bv)\rangle_s{\langle I,I\rangle_s}
\;\ge\;
\langle I,A(\bv)\rangle_s^2
\]
since \(\langle I,I\rangle_s>0\).  But
\[
\langle I,A(\bv)\rangle_s
=\Tr \bigl(I\,X^s(\by)\,A(\bv)\,X^{1-s}(\by)\bigr)
=\Tr \bigl(A(\bv)\,X(\by)\bigr),
\]
and
\[
\langle I,I\rangle_s
=\Tr \bigl(I\,X^s(\by)\,I\,X^{1-s}(\by)\bigr)
=\Tr X(\by)
=1.
\]
Hence
\begin{equation}\label{eq5}
	\langle A(\bv),A(\bv)\rangle_s = \Tr \bigl(A(\bv)\,X^s(\by)\,A(\bv)\,X^{1-s}(\by)\bigr)
	\;\ge\;
	\bigl[\Tr\bigl(A(\bv)\,X(\by)\bigr)\bigr]^2.
\end{equation}
Integrating over \(s\in[0,1]\) gives
\[
\int_0^1
\Tr \bigl(A(\bv)X^s(\by)A(\bv)X^{1-s}(\by)\bigr)\,ds
\;\ge\;
\int_0^1\bigl[\Tr(A(\bv)\,X(\by))\bigr]^2\,ds
=\bigl[\Tr(A(\bv)\,X(\by))\bigr]^2.
\]
Therefore
\[
\bv^T\nabla^2 f(\by)\,\bv
\;=\;
\int_0^1\Tr(A(\bv)X^s(\by)A(\bv)X^{1-s}(\by))\,ds
\;-\;\bigl[\Tr(A(\bv)\,X(\by))\bigr]^2
\;\ge\;0,
\]
and so \(\nabla^2 f(\by)\) is positive semidefinite.\end{proof}

\begin{lemma}\label{coercive}
Function \(f\) is coercive (i.e.\ \(\lim_{\|y\|\to\infty}f(y)=+\infty\)) 
if and only if \(\bb \in \Int \mathscr M\).
\end{lemma}

\begin{proof}
{	Let $h$ denote the support function of $\mathscr M$ defined in \eqref{eq:support}, and define $g(\bv):=h(\bv)-\bb^T\bv$, a continuous, positively
	homogeneous function on $\R^m$.
	
	Fix $\bv\in\R^m$ and let $\mu_1\ge\cdots\ge\mu_n$ be the eigenvalues of $A(\bv)$, so that
	$\mu_1=h(\bv)$. For every $t\ge0$,
	\[
	e^{t\mu_1}\;\le\;\Tr\exp(tA(\bv))=\sum_{k=1}^n e^{t\mu_k}\;\le\;n\,e^{t\mu_1},
	\]
	and taking logarithms and subtracting $t\,\bb^T\bv$ yields, along the ray $\by=t\bv$,
	\begin{equation}\label{eq:raybound}
		t\,g(\bv)\;\le\;f(t\bv)\;\le\;t\,g(\bv)+\log n,\qquad t\ge0.
	\end{equation}
	
	Since $\mathscr M$ is convex and compact with support function $h$, one has
	$\bb\in\mathscr M$ iff $\bb^T\bv\le h(\bv)$ for all $\bv$, and
	\begin{equation}\label{eq:intchar}
		\bb\in\Int\mathscr M\quad\Longleftrightarrow\quad g(\bv)>0\ \text{ for all }\bv\neq0.
	\end{equation}
	Indeed, if $\bb\in\Int\mathscr M$ then $\bb+\varepsilon\bv/\|\bv\|\in\mathscr M$ for some
	$\varepsilon>0$, so $h(\bv)\ge\bb^T\bv+\varepsilon\|\bv\|>\bb^T\bv$. Conversely, if $g>0$ on the
	unit sphere $\mathscr V:=\{\bv\in\R^m:\|\bv\|=1\}$, then by continuity and compactness
	$$\gamma:=\min_{\bv\in\mathscr V}g(\bv)>0$$ and for any $\bw$ with $\|\bw\|\le\gamma$ and any
	$\bv\in\mathscr V$,
	\[
	(\bb+\bw)^T\bv=\bb^T\bv+\bw^T\bv\le(h(\bv)-\gamma)+\|\bw\|\le h(\bv),
	\]
	whence $\bb+\bw\in\mathscr M$, and therefore $\bb\in\Int\mathscr M$. By homogeneity it suffices
	to evaluate $g$ on $\mathscr V$.
	
	\smallskip{\sl Sufficiency.}
	If $\bb\in\Int\mathscr M$, then $\gamma>0$ by \eqref{eq:intchar}.
	For $\by\neq0$, putting $t=\|\by\|$ and $\bv=\by/\|\by\|\in\mathscr V$, the left inequality in
	\eqref{eq:raybound} gives
	\[
	f(\by)\ge\|\by\|\,g(\bv)\ge\gamma\,\|\by\|\;\xrightarrow[\|\by\|\to\infty]{}\;+\infty,
	\]
	so $f$ is coercive.
	
	\smallskip\noindent{\sl Necessity.}
	If $\bb\notin\Int\mathscr M$, then by \eqref{eq:intchar} there is $\bv_0\in\mathscr V$ with
	$g(\bv_0)\le0$. The right inequality in \eqref{eq:raybound} along $\by=t\bv_0$ gives
	\[
	f(t\bv_0)\le t\,g(\bv_0)+\log n\le\log n\qquad(t\ge0),
	\]
	so $f$ remains bounded above on the unbounded ray $\{t\bv_0:t\ge0\}$ and cannot tend to
	$+\infty$ as $\|\by\|\to\infty$; hence $f$ is not coercive. }
	\end{proof}

\begin{lemma}\label{residualdist}
	For every $\bb\in\R^m$,
	\[
	\inf_{\by\in\R^m}\|\nabla f(\by)\|_2=\mathrm{dist}(\bb,\mathscr M).
	\]
\end{lemma}

\begin{proof}
	By Lemma~\ref{gh}, $\nabla f(\by)=\mathcal A(X(\by))-\bb$ with $X(\by)=\exp_1 A(\by)\in\S^n_1$, so
	$\mathcal A(X(\by))\in\mathscr M$ for every $\by$ and
	\[
	\|\nabla f(\by)\|_2=\|\mathcal A(X(\by))-\bb\|_2\ \ge\ \mathrm{dist}(\bb,\mathscr M).
	\]
	Taking the infimum over $\by$ gives	$\inf_{\by\in\R^m}\|\nabla f(\by)\|_2 \geq \mathrm{dist}(\bb,\mathscr M).$
	
	For the other direction, let $\bb^\star$ be the nearest point of the
	compact convex set $\mathscr M$ to $\bb$, so that $\|\bb^\star-\bb\|_2=\mathrm{dist}(\bb,\mathscr M)$,
	and let $\bb^\circ:=\mathcal A(\tfrac1n I)\in\Int\mathscr M$ be the center of $\mathscr M$, the image of
	the maximally mixed state. For $\epsilon\in(0,1]$ set
	\[
	\bb^\epsilon:=(1-\epsilon)\,\bb^\star+\epsilon\,\bb^\circ .
	\]
	Since $\bb^\star\in\mathscr M$ and $\bb^\circ\in\Int\mathscr M$, it holds
	$\bb^\epsilon\in\Int\mathscr M$. By Lemma~\ref{coercive} the function
	$\by\mapsto\log\Tr\exp A(\by)-(\bb^\epsilon)^T\by$ is then coercive and continuous, hence attains its
	minimum at some $\by^\epsilon\in\R^m$; being convex and differentiable, its gradient vanishes there, i.e.
	$A(X(\by^\epsilon))=\bb^\epsilon$.
	Therefore $\nabla f(\by^\epsilon)=\mathcal A(X(\by^\epsilon))-\bb=\bb^\epsilon-\bb$, and
	\[
	\inf_{\by}\|\nabla f(\by)\|_2\ \le\ \|\bb^\epsilon-\bb\|_2
	=\|(1-\epsilon)\bb^\star+\epsilon\,\bb ^\circ-\bb\|_2
	\ \xrightarrow[\ \epsilon\to0^+\ ]{}\ \|\bb^\star-\bb\|_2=\mathrm{dist}(\bb,\mathscr M).
	\]
	Combining the two inequalities proves the identity.
\end{proof}	

\begin{lemma}\label{unique}
If $\bb \in \Int \mathscr M$, function $f$ has a unique global minimizer
\[
{\by^*} = \arg \min_{\by \in \R^m} f(\by).
\]
\end{lemma}

\begin{proof}

{Recall expression \eqref{eq4} from the proof of Lemma \ref{gh}, let
\[
V(\by, \bv) := \int_0^1
\Tr(A(\bv)\,X^s(\by)\,A(\bv)\,X^{1-s}(\by))\,ds
-
\bigl[\Tr\bigl(A(\bv)\,X(\by)\bigr)\bigr]^2
\]
and let us show that $V(\by, \bv)$ is zero if and only if $A(\bv)$ is a multiple of the identity matrix.

Assume $A(\bv) = cI$ for some scalar $c \in \mathbb{R}$. We substitute this into the expression for $V(\by, \bv)$.
The first term becomes:
\[
\int_0^1 \Tr\bigl((cI)\,X^s(\by)(cI)X^{1-s}(\by)\bigr)\,ds = \int_0^1 \Tr\bigl(c^2 X^s(\by) X^{1-s}(\by)\bigr)\,ds = \int_0^1 c^2 \Tr X(\by)\,ds = c^2.
\]
The second term becomes:
\[
\bigl[\Tr\bigl((cI) X(\by)\bigr)\bigr]^2 = \bigl[c \Tr(X(\by))\bigr]^2 = c^2.
\]
Therefore, $V(\by, \bv) = 0$.

The converse statement relies on the Cauchy-Schwarz inequality \eqref{eq5} obtained from the inner product $\langle .,. \rangle_s$ defined in the proof of Lemma \ref{gh}.
This inequality shows that the integrand in the expression for $V(\by, \bv)$ is non-negative for all $s \in [0,1]$. If $V(\by, \bv) = 0$ it holds
\[
\int_0^1 \left( \Tr(A(\bv)X^s(\by) A(\bv) X^{1-s}(\by))\right) ds = 
\bigl[\Tr\bigl(A(\bv)\,X(\by)\bigr)\bigr]^2.
\]
Since the integrand is a continuous and non-negative function of $s$, its integral can only be zero if the integrand is identically zero for all $s \in [0,1]$. Therefore:
\[
\Tr( A(\bv)X^s(\by) A(\bv) X^{1-s}(\by)) = \bigl[\Tr(A(\bv)X(\by))\bigr]^2 \quad \text{for all } s \in [0,1].
\]
This means the Cauchy-Schwarz inequality \eqref{eq5} must hold with equality for all $s$. Equality in the Cauchy-Schwarz inequality holds if and only if the two matrices $A(\bv)$ and $I$ are linearly dependent, meaning $A(\bv)$ must be a multiple of the identity matrix. 

Under Assumption 1, the matrices $A_1, \ldots, A_m$ cannot span the identity matrix, and hence $A(\bv) = cI$ implies that $c=0$, and thus $A(\bv) = \mathbf{0}$. The linear independence of the matrices then forces $\bv = \mathbf{0}$.
 
 Therefore, $\bv^T\nabla^2 f(\by)\,\bv > 0$ for all $\bv \neq \mathbf{0}$, proving that $f$ is strictly convex. A strictly convex function has at most one minimizer. Since $f$ is also coercive when $\bb \in \Int \mathscr M$ (by Lemma \ref{coercive}), it is guaranteed to have a unique global minimizer.
}
\end{proof}

The above results suggest that minimizing $f$ solves the moment body membership oracle {when $\bb \in \Int\mathscr M$.} Indeed, $f$ has a unique global minimizer $\by^*$ at which the gradient of $f$ vanishes, i.e. $\mathcal A(X(\by^*))=\bb$. Therefore the inclusion $\bb\in \mathscr M$ is certified by the matrix $X(\by^*)=\exp_1 A(\by^*) \in \S^n_1$.

{
\begin{theorem}[maximum entropy duality]\label{primaldual}
	If $\bb\in\Int\mathscr M$, the convex unconstrained minimization
	\[\boxed{
	d^\star:=\min_{\by\in\R^m} f(\by),\qquad
	f(\by)=\log\trace\exp A(\by)-\bb^T\by,
	}\]
	and the concave maximum entropy problem
	\[\boxed{
	p^\star:=\max_{X\in\S^n_1}\,S(X)\quad\text{s.t.}\quad \mathcal A(X)=\bb,
	\qquad S(X):=-\trace(X\log X),
	}\]
	are related by strong duality, i.e.  $p^\star=d^\star$,
	and both extrema are attained. The dual minimizer $\by^\star$ is unique, and
	the unique primal maximizer is
	\[
	\boxed{\,X^\star=\exp_1 A(\by^\star)\;\succ\;0\,,}
	\]
	so in particular $\mathcal A(X^\star)=\bb$.
\end{theorem}

\begin{proof}
	Fix $Z\in\S^n$ and let $\phi(X):=S(X)+\langle X,Z\rangle$
	on $\S^n_1$. With $X':=\exp_1 Z\succ0$ one has $\log X'=Z-(\log\trace\exp Z)\,I$, so, using
	$\nabla_X\trace(X\log X)=\log X+I$,
	\[
	\nabla\phi(X')=Z-\log X'-I=(\log\trace\exp Z-1)\,I .
	\]
	Being a multiple of $I$, $\nabla\phi(X')$ is orthogonal to $X-X'$ for every $X\in\S^n_1$
	(both have unit trace). Since $S$ and hence $\phi$ is strictly concave 	\cite[eq.~(1.11)]{carlen10} and differentiable at $X'$, the gradient
	inequality therefore gives, for all $X\in\S^n_1$,
	\[
	\phi(X)\;\le\;\phi(X')+\langle\nabla\phi(X'),X-X'\rangle=\phi(X'),
	\]
	with equality iff $X=X'$.
 Since $\phi(X')=\trace(X'(Z-\log X'))=\log\trace\exp Z$,
	\begin{equation}\label{eq:varineq}
		S(X)+\langle X,Z\rangle\;\le\;\log\trace\exp Z,
		\qquad\text{equality iff }X=\exp_1 Z .
	\end{equation}
	
	For feasible $X$  and any $\by\in\R^m$, set
	$Z=A(\by)$ in \eqref{eq:varineq} and use $\langle X,A(\by)\rangle=\by^T\mathcal A(X)=\bb^T\by$:
	\[
	S(X)\;\le\;\log\trace\exp A(\by)-\bb^T\by=f(\by),
	\]
	hence $p^\star\le d^\star$.
	
	By Lemmas~\ref{coercive} and~\ref{unique}, the hypothesis
	$\bb\in\Int\mathscr M$ makes $f$ coercive and strictly convex, so $f$ has a unique minimizer
	$\by^\star$, with $\nabla f(\by^\star)=\mathcal A(\exp_1 A(\by^\star))-\bb=0$. Let
	$X^\star:=\exp_1 A(\by^\star)$; then $X^\star\succ0$ and $\mathcal A(X^\star)=\bb$, so $X^\star$ is
	feasible. The equality case of \eqref{eq:varineq} with $Z=A(\by^\star)$ gives
	\[
	S(X^\star)=\log\trace\exp A(\by^\star)-\bb^T\by^\star=f(\by^\star)=d^\star ,
	\]
	and feasibility gives $p^\star\ge S(X^\star)=d^\star$. With $p^\star\le d^\star$ this forces
	$p^\star=d^\star=S(X^\star)$: the matrix $X^\star$ attains the primal maximum, and it is the only
	maximizer since $S$ is strictly concave on $\S^n_1$.
\end{proof}
}

The following result is well-known in quantum information theory, see e.g. \cite[Theorem 2]{jv23}, where it is attributed to \cite{w63}.

\begin{lemma}\label{diffeo}
The map $\by \mapsto {\mathcal A}(\exp_1  A(\by))$ is a smooth diffeomorphism between $\R^m$ and  $\Int \mathscr M$.
\end{lemma}

\begin{proof}
	Let $g(\by):=\log\trace\exp A(\by)$. By Lemma~\ref{gh}, $g$ is smooth on $\R^m$ and,
	since $g$ and $f$ differ only by the linear term $\bb^T\by$,
	\[
	\Phi(\by):=\nabla g(\by)=\mathcal A(X(\by)),\qquad
	\nabla^2 g(\by)=\nabla^2 f(\by),\qquad
	X(\by)=\exp_1 A(\by)\succ0 .
	\]
	Under Assumption~\ref{injective} the Hessian is positive definite at every point: the
	strict-convexity argument in the proof of Lemma~\ref{unique} gives
	$\bv^T\nabla^2 f(\by)\bv>0$ for all $\bv\neq0$ and all $\by$, using only the linear
	independence of $I_n,A_1,\dots,A_m$ and not the position of $\bb$. Thus $g$ is smooth and
	strictly convex on $\R^m$.
	
	\smallskip\noindent{(i) $\Phi$ is an injective local diffeomorphism.}
	For every $\by$, $\nabla\Phi(\by)=\nabla^2 g(\by)\succ0$ is invertible, so by the inverse
	function theorem $\Phi$ is a local diffeomorphism; in particular it is an open map. Strict
	convexity of $g$ yields strict monotonicity of its gradient,
	\[
	(\Phi(\by_1)-\Phi(\by_2))^{T}(\by_1-\by_2)>0,\qquad\by_1\neq\by_2,
	\]
	so $\Phi$ is injective. An injective local diffeomorphism is a diffeomorphism onto its
	image, and that image $\Phi(\R^m)$ is open.
	
	\smallskip\noindent\textbf{(ii) $\Phi(\R^m)\subseteq\Int\mathscr M$.}
	For every $\by$ one has $X(\by)\in\S^n_1$, hence
	$\Phi(\by)=\mathcal A(X(\by))\in\mathcal A(\S^n_1)=\mathscr M$. Since $\Phi(\R^m)$ is open
	by (i), it is an open subset of $\mathscr M$ and therefore contained in $\Int\mathscr M$.
	
	\smallskip\noindent\textbf{(iii) $\Phi(\R^m)\supseteq\Int\mathscr M$.}
	Fix $\bx\in\Int\mathscr M$ and set
	\[
	f_\bx(\by):=g(\by)-\bx^T\by=\log\trace\exp A(\by)-\bx^T\by ,
	\]
	which is precisely the function of Lemma~\ref{coercive} with target $\bx$ in place of
	$\bb$. As $\bx\in\Int\mathscr M$, Lemma~\ref{coercive} shows $f_\bx$ is coercive, and it is
	strictly convex because it differs from $g$ only by a linear term. A coercive, smooth,
	strictly convex function attains its minimum at a unique point $\by$, where
	$\nabla f_\bx(\by)=\Phi(\by)-\bx=0$, i.e.\ $\Phi(\by)=\bx$.
	
	\smallskip
	By (ii) and (iii), $\Phi(\R^m)=\Int\mathscr M$; combined with (i), $\Phi$ is a
	diffeomorphism from $\R^m$ onto $\Int\mathscr M$.
\end{proof}

\begin{figure}[h]
	\centering
	
	\begin{minipage}{0.45\textwidth}
		\centering
		\includegraphics[width=\linewidth]{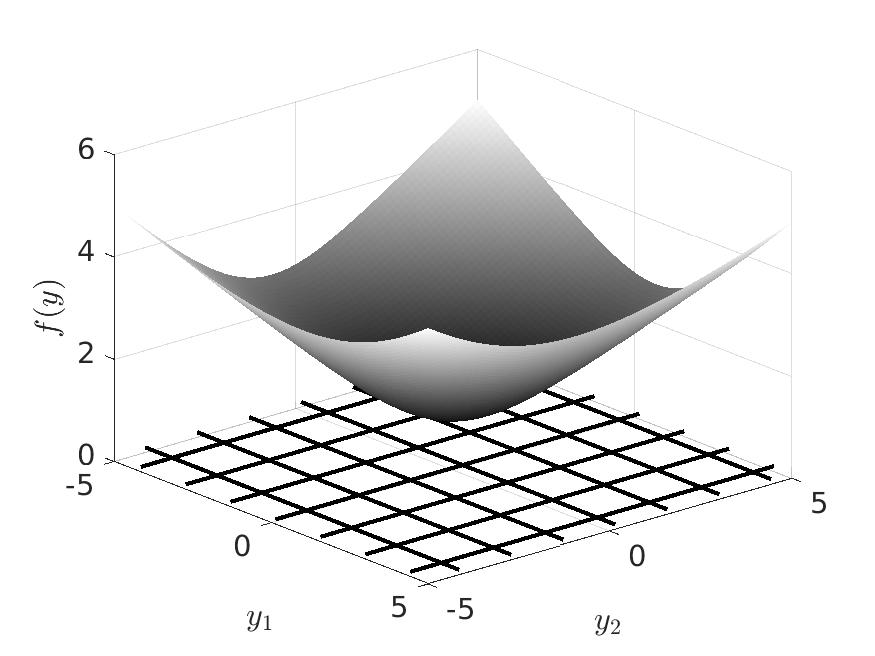}
		\caption{Dual function graph (gray) and regular grid underneath (black).}
		\label{fig:dualfunction}
	\end{minipage}\hspace{1em}
	\begin{minipage}{0.45\textwidth}
		\centering
		\includegraphics[width=\linewidth]{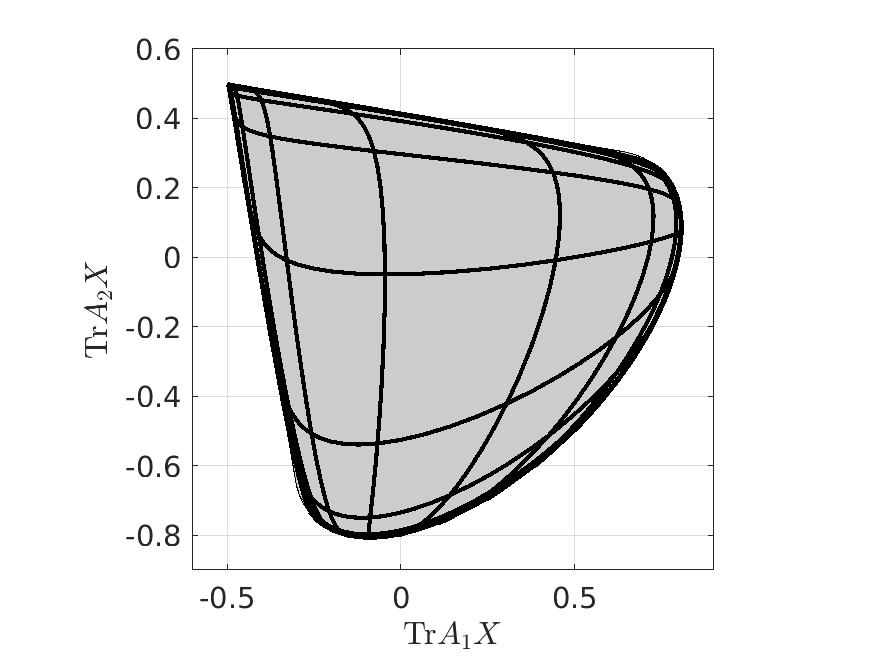}
		\caption{Moment body (gray) and grid image through the gradient map (black).}
		\label{fig:momentbodygrid2d}
	\end{minipage}
	
\end{figure}

\begin{example}\label{ex:2d}
For the matrices of Example \ref{ex:momentbody2d},
the graph of function $f(y)$ is represented on Figure \ref{fig:dualfunction}, together with a regular grid $\mathscr Y$ (black lines underneath).
The moment body $\mathscr M$ is represented on Figure \ref{fig:momentbodygrid2d}, together with the image of the grid through the gradient map $\mathcal A(\exp_1 A(\mathscr Y))$.

\end{example}

\begin{remark}[Commutative specialization and classical counterparts]
	Theorem~\ref{primaldual} and Lemma~\ref{diffeo} are the noncommutative (matrix) versions of
	classical results on entropy-regularized problems over probability vectors and densities, and it
	is instructive to record the commutative specialization. When the data matrices $A_1,\dots,A_m$
	commute they are simultaneously diagonalizable, $X(\by)=\exp_1 A(\by)$ is diagonal in the common
	eigenbasis, and all spectral objects above reduce to their scalar analogues: the spectraplex
	$\S^n_1$ becomes the probability simplex $\Delta=\{\bp\ge0:\sum_k p_k=1\}$, the von Neumann
	entropy $S(X)=-\trace(X\log X)$ becomes the Shannon entropy $H(\bp)=-\sum_k p_k\log p_k$, and the
	log-partition function $\log\trace\exp A(\by)$ becomes the log-sum-exp
	$\log\sum_k\exp(\sum_i y_i a_{ik})$, where $a_{ik}$ is the $k$-th eigenvalue of $A_i$. In this
	regime Theorem~\ref{primaldual} is exactly the  maximum-entropy principle:
	maximizing $H(\bp)$ subject to the moment constraints $\sum_k p_k a_{ik}=b_i$ is dual to
	minimizing the smooth convex log-sum-exp, the unique optimizer being the Gibbs (exponential-family)
	law $p^\star_k = \exp(\sum_i y^\star_i a_{ik}) / \sum_j  \exp(\sum_i y^\star_i a_{ij})$. This is the standard conjugacy between the
	negative Shannon entropy and the log-sum-exp \cite{BoydVandenberghe04}, which underlies
	the maximum-entropy method \cite{ct06}. The continuous-density version, in which $\Delta$ is replaced by
	probability densities on a set and $H$ by differential entropy, is the classical maximum-entropy
	moment problem; its strong duality and the exponential-family form of the optimal density are due
	to \cite{bl91}.

	Likewise, Lemma~\ref{diffeo} is the matrix form of a cornerstone of the theory of exponential
	families. The gradient of the log-partition function---the forward or moment map sending the
	canonical (natural) parameters $\by$ to the expectation (mean) parameters $\mathcal A(X(\by))$---is
	a diffeomorphism from the natural parameter space onto the relative interior of the moment
	(marginal) polytope, its inverse being the Legendre transform conjugate to the entropy; see
	\cite[Ch.~3]{wj08} and, from the information-geometric
	viewpoint, the dually flat structure of \cite[Ch.~3]{an00}. The curvature analysis of the next
	sections specializes just as transparently: in the commutative case the
	Hessian of Lemma~\ref{gh} collapses (the integral over $s$ disappears) to the ordinary covariance
	matrix of the moment vector under the Gibbs law $\bp$,
	that is, the Fisher information matrix of the exponential family. The best-known large-scale instance of this commutative picture
	is entropically regularized optimal transport, where membership in the relative interior of a
	transportation polytope is decided by minimizing precisely such a log-sum-exp dual through Sinkhorn
	iterations \cite{cuturi13,pc19}; the Shannon- and von Neumann-entropy regularizations of general
	linear and semidefinite programs in \cite{l23,cl25b} are the noncommutative extension of the same
	idea, of which our trace-one membership oracle is a special case.
\end{remark}

\section{Geometric analysis}\label{curvature}

Since $f$ is smooth and convex, its minimization can be achieved with standard optimization algorithms.
The performance of these algorithms depends on the geometry of $f$, and especially its curvature.
We say that $f$ is $\alpha$-strongly convex and $\lambda$-smooth whenever
\[
0 \prec \alpha I_m \preceq \nabla^2f(\by)\preceq \lambda I_m
\quad\forall\,\by \in \R^m.
\quad
\]
The constant $\alpha>0$ is called the strong convexity modulus, and the constant $\lambda>0$ is called the smoothness constant (or Lipschitz constant of $\nabla f$).
The condition number
\[
\kappa=\frac{\lambda}{\alpha}
\] 
governs the convergence rates of standard first‐ and second‐order methods.
With fixed step‐size \(1/\lambda\), gradient descent 
\(\;\by_{k+1}=\by_k-{1/\lambda}\nabla f(\by_k)\)
satisfies the linear rate
$
f(\by_k)-f(\by^*)
\;\le\;(1-\frac{1}{\kappa})^k (f(\by_0)-f(\by^*)),
$
so that reaching \(\varepsilon\)-accuracy to the minimum $\by^*$ requires
\(\,O(\kappa\log\varepsilon^{-1})\) iterations, see 
e.g. \cite[Sec.\,9.3.1]{BoydVandenberghe04}. 
Nesterov's accelerated scheme achieves the optimal first-order complexity
\(O(\sqrt{\kappa}\log(1/\varepsilon))\) by combining momentum with gradient
steps; see \cite[Ch.~2, Thm.~2.2.2]{Nesterov18}.
Near the optimum, Newton’s method
\(\;\by_{k+1}=\by_k - [\nabla^2f(\by_k)]^{-1}\nabla f(\by_k)\)
converges quadratically, but its {region of attraction} and the quality
of each step depend on \(\kappa\).  Ill‐conditioned Hessians can force
small steps or necessitate line‐search/globalization strategies, whose
complexity again scales with \(\kappa\), see   \cite[Sec.\,3.5]{NocedalWright06}.
Interior‐point methods exhibit polynomial‐time complexity bounds that
depend on the barrier Hessian’s conditioning \(\kappa\)
(see \cite[Chap.\,5]{Wright94IPM}), and quasi‐Newton updates
(e.g.\ BFGS) achieve superlinear convergence only when \(\kappa\)
is moderate. 
{As shown in \cite[Thm.~2 and Cor.~1]{JinJiangMokhtari2025}, 
for exact line-search BFGS, the objective values converge globally linearly at the rate $(1-1/(3\kappa))^k$.
}

To solve the moment body membership problem for moderate size problems, we propose to use L-BFGS, a standard quasi-Newton algorithm constructing an approximation of the Hessian using a limited number of evaluations of the gradient. It can be interpreted as a discretization of a variable-metric generalization of the Newton flow 	where the true inverse Hessian is replaced by a time-varying symmetric positive-definite matrix.

In the context of semidefinite optimization, the idea of formulating and solving with BFGS a dual smooth problem was already explored in \cite{m04} for the semidefinite least-squares problem, consisting of projecting a given symmetric matrix onto a given spectrahedral shadow. It was later on used to solve polynomial SOS problems \cite{hm11}.

\subsection{Smoothness}

{
In this subsection we bound the curvature of $f$ from above. The relevant data is the Gram matrix
of the traceless parts of the $A_i$. Defining
\[
A_i^0:=A_i-\frac{\Tr A_i}{n}\,I_n,
\]
we define the centered Gram matrix
\begin{equation}\label{gram0}
	G_0:=\bigl[\Tr(A_i^0A_j^0)\bigr]_{i,j=1}^m,
\end{equation}
with eigenvalues $\lambda_{\max}(G_0)\ge\cdots\ge\lambda_{\min}(G_0)$. Since
$A(\bv)-\tfrac{\Tr A(\bv)}{n}I_n=\sum_iv_iA_i^0$ is the orthogonal projection of $A(\bv)$ onto the
trace-zero subspace of $\S^n$, one has, for every $\bv\in\R^m$,
\begin{equation}\label{g0id}
	\bv^TG_0\bv
	=\bigl\|A(\bv)-\tfrac{\Tr A(\bv)}{n}I_n\bigr\|_F^2
	=\Tr(A(\bv)^2)-\tfrac1n\bigl(\Tr A(\bv)\bigr)^2\;\ge\;0 .
\end{equation}
Under Assumption~\ref{injective}, $G_0$ is positive definite: if $\bv^TG_0\bv=0$ then $A(\bv)$ is a
scalar multiple of $I_n$, and linear independence of $I_n,A_1,\dots,A_m$ forces $\bv=0$. Hence
$\lambda_{\min}(G_0)>0$.

\begin{lemma}[Smoothness]\label{smoothness}
	For every $\by\in\R^m$ one has $$\nabla^2f(\by)\preceq\tfrac12 G_0.$$ In other words, $f$ is
	$\tfrac12\lambda_{\max}(G_0)$-smooth.
\end{lemma}

\begin{proof}
	Fix $\by,\bv\in\R^m$ and write, for brevity, $X:=X(\by)\succ0$ and $A:=A(\bv)\in\S^n$.
	
	By the directional Hessian \eqref{eq4},
	\[
	\bv^T\nabla^2f(\by)\,\bv=\int_0^1 a(s)\,ds-\bigl[\Tr(AX)\bigr]^2,
	\qquad a(s):=\Tr(A\,X^s\,A\,X^{1-s}).
	\]
	Diagonalize $X=\sum_{i=1}^n\lambda_i\,\bq_i\bq_i^T$ with $\lambda_i>0$ and orthonormal $\bq_i$. Since
	$A=A^T$, for every $s\in[0,1]$,
	\[
	a(s)=\sum_{i,j=1}^n(\bq_i^TA\bq_j)^2\,\lambda_i^{\,s}\lambda_j^{\,1-s}
	=\sum_{i,j=1}^n(\bq_i^TA\bq_j)^2\,\lambda_j\,e^{\,s\log(\lambda_i/\lambda_j)} .
	\]
	Each summand is a nonnegative multiple of a convex exponential of $s$, so $a$ is convex on $[0,1]$,
	and $a(0)=a(1)=\Tr(A^2X)$. A convex function lies below the chord through its endpoints, here the
	constant $\Tr(A^2X)$, whence $\int_0^1 a(s)\,ds\le\Tr(A^2X)$ and
	\begin{equation}\label{upperbound}
		\bv^T\nabla^2f(\by)\,\bv\;\le\;\Tr(A^2X)-\bigl[\Tr(AX)\bigr]^2 .
	\end{equation}
	
	For every $c\in\R$, using $\Tr X=1$,
	\[
	\Tr((A-cI_n)^2X)=\Tr(A^2X)-2c\,\Tr(AX)+c^2,
	\qquad
	\bigl[\Tr((A-cI_n)X)\bigr]^2=\bigl[\Tr(AX)-c\bigr]^2 ,
	\]
	and subtracting the second identity from the first gives, for every $c$,
	\[
	\Tr(A^2X)-\bigl[\Tr(AX)\bigr]^2
	=\Tr((A-cI_n)^2X)-\bigl[\Tr(AX)-c\bigr]^2
	\;\le\;\Tr((A-cI_n)^2X),
	\]
	discarding the nonnegative square. For any $M\in\S^n$ with $M\succeq0$ one has
	$\lambda_{\max}(M)I_n-M\succeq0$, hence $0\le\Tr((\lambda_{\max}(M)I_n-M)X)=\lambda_{\max}(M)-\Tr(MX)$,
	i.e.\ $\Tr(MX)\le\lambda_{\max}(M)$. Applying this to $M=(A-cI_n)^2\succeq0$, whose eigenvalues are
	$(\lambda_i(A)-c)^2$,
	\[
	\Tr(A^2X)-\bigl[\Tr(AX)\bigr]^2\;\le\;\lambda_{\max}((A-cI_n)^2)=\max_{i}(\lambda_i(A)-c)^2 .
	\]
	Choosing $c=\tfrac12(\lambda_{\max}(A)+\lambda_{\min}(A))$ makes the right-hand side smallest, equal
	to $\tfrac14(\lambda_{\max}(A)-\lambda_{\min}(A))^2$, so by \eqref{upperbound},
	\begin{equation}\label{spectralbound}
		\bv^T\nabla^2f(\by)\,\bv\;\le\;\tfrac14\bigl(\lambda_{\max}(A)-\lambda_{\min}(A)\bigr)^2 .
	\end{equation}
	
	The spread $\lambda_{\max}(A)-\lambda_{\min}(A)$ is unchanged when $A$ is shifted by a multiple of
	$I_n$, hence equals the spread of the traceless part $A^0(\bv)=A-\tfrac{\Tr A}{n}I_n$. Let
	$\nu_1\ge\cdots\ge\nu_n$ be the eigenvalues of $A^0(\bv)$, so that
	$\sum_k\nu_k^2=\Tr(A^0(\bv)^2)=\bv^TG_0\bv$ by \eqref{g0id}. Then $0\le(\nu_1+\nu_n)^2$ gives
	\[
	\bigl(\lambda_{\max}(A)-\lambda_{\min}(A)\bigr)^2=(\nu_1-\nu_n)^2
	\le2(\nu_1^2+\nu_n^2)\le2\sum_k\nu_k^2=2\,\bv^TG_0\bv .
	\]
	Combining with \eqref{spectralbound}, $\bv^T\nabla^2f(\by)\,\bv\le\tfrac12\,\bv^TG_0\bv$ for every
	$\bv\in\R^m$, which is the operator inequality $\nabla^2f(\by)\preceq\tfrac12 G_0$. In particular
	$\bv^T\nabla^2f(\by)\,\bv\le\tfrac12\lambda_{\max}(G_0)\|\bv\|_2^2$, so $f$ is
	$\tfrac12\lambda_{\max}(G_0)$-smooth.
\end{proof}

\begin{remark}
	The bound \eqref{spectralbound},
	\[
	\Tr(A(\bv)^2\,X(\by)) - [\Tr(A(\bv)X(\by))]^2
	\;\le\;\tfrac14\bigl(\lambda_{\max}(A(\bv))-\lambda_{\min}(A(\bv))\bigr)^2,
	\]
	proved here algebraically, is a particular case of a more general result in non-commutative
	probability theory, see \cite[Thm.~2]{BD2000}.
\end{remark}
}

\subsection{Strong convexity}

{
We now bound the curvature of $f$ from below, in terms of the same centered Gram matrix $G_0$
introduced in \eqref{gram0}. In contrast with the uniform upper bound $\tfrac12 G_0$ of
Lemma~\ref{smoothness}, the lower bound carries the smallest eigenvalue of the density matrix
$X(\by)$, which degrades as $\by$ moves away from the origin.

\begin{lemma}[Pointwise curvature]\label{curvature}
	For every $\by\in\R^m$, with $X(\by)=\exp_1A(\by)\succ0$,
	\[
	\nabla^2f(\by)\;\succeq\;\lambda_{\min}(X(\by))\,G_0.
	\]
\end{lemma}

\begin{proof}
	Fix $\by\in\R^m$ and $\bv\in\R^m$. Diagonalize $X(\by)=\sum_{k=1}^n p_k\,\bq_k\bq_k^T$, with eigenvalues
	$p_k=\lambda_k(X(\by))>0$ satisfying $\sum_kp_k=1$, orthonormal eigenvectors $\bq_k$, and write
	$\rho:=\lambda_{\min}(X(\by))=\min_kp_k$. Introduce the entries of $A(\bv)$ in this eigenbasis,
	\[
	a_{ij}:=\bq_i^T A(\bv)\,\bq_j\qquad i,j=1,\dots,n .
	\]
	Expanding the directional Hessian \eqref{eq4} in this eigenbasis exactly as in the proof of
	Lemma~\ref{smoothness} gives
	\[
	\int_0^1\Tr(A(\bv)\,X(\by)^s\,A(\bv)\,X(\by)^{1-s})\,ds
	=\sum_{i,j=1}^n a_{ij}^2\,\ell(p_i,p_j),\qquad \ell(p_i,p_j):=\int_0^1 p_i^{\,s}p_j^{\,1-s}\,ds,
	\]
	while $\Tr(A(\bv)\,X(\by))=\sum_{k=1}^n p_k\,a_{kk}$. Using $\ell(p_k,p_k)=p_k$ to split off the
	diagonal, \eqref{eq4} reads
	\[
	\bv^T\nabla^2f(\by)\,\bv
	=\underbrace{\sum_{k} p_k\,a_{kk}^2-\Bigl(\sum_{k} p_k\,a_{kk}\Bigr)^{\!2}}_{T_1}
	\;+\;\underbrace{\sum_{i\ne j} a_{ij}^2\,\ell(p_i,p_j)}_{T_2}.
	\]
	For $s\in[0,1]$ the weighted geometric mean dominates the smaller argument,
	$p_i^{\,s}p_j^{\,1-s}\ge\min(p_i,p_j)\ge\rho$, hence $\ell(p_i,p_j)\ge\rho$ and
	$T_2\ge\rho\sum_{i\ne j} a_{ij}^2$. Setting $\bar m:=\sum_k p_k\,a_{kk}$, the term
	$T_1=\sum_k p_k\,(a_{kk}-\bar m)^2$ is a weighted variance, so
	\[
	T_1\;\ge\;\rho\sum_k (a_{kk}-\bar m)^2
	\;\ge\;\rho\,\min_{c\in\R}\sum_k (a_{kk}-c)^2
	=\rho\Bigl(\sum_k a_{kk}^2-\tfrac1n\bigl(\textstyle\sum_k a_{kk}\bigr)^2\Bigr),
	\]
	the inner minimum attained at $c=\tfrac1n\sum_k a_{kk}$. Adding the two bounds,
	\[
	\bv^T\nabla^2f(\by)\,\bv\;\ge\;\rho\Bigl(\sum_{i,j} a_{ij}^2-\tfrac1n\bigl(\textstyle\sum_k a_{kk}\bigr)^2\Bigr).
	\]
	By orthogonal invariance of the Frobenius norm, $\sum_{i,j}a_{ij}^2=\|Q^TA(\bv)Q\|_F^2=\Tr(A(\bv)^2)$
	with $Q:=[\bq_1\cdots \bq_n]$, and $\sum_k a_{kk}=\Tr(Q^TA(\bv)Q)=\Tr A(\bv)$. Hence
	\[
	\bv^T\nabla^2f(\by)\,\bv
	\;\ge\;\rho\Bigl(\Tr(A(\bv)^2)-\tfrac1n(\Tr A(\bv))^2\Bigr)
	=\lambda_{\min}(X(\by))\,\bv^TG_0\bv,
	\]
	the last equality by \eqref{g0id}. As $\bv\in\R^m$ was arbitrary and both $\nabla^2f(\by)$ and
	$\lambda_{\min}(X(\by))\,G_0$ are symmetric, this is the asserted inequality.
\end{proof}

\begin{remark}
	At $\by=0$ one has $X(0)=\tfrac1n I_n$, and the curvature bound holds with equality:
	$\nabla^2f(0)=\tfrac1n G_0$. Lemma~\ref{curvature} is therefore sharp, and together with
	Lemma~\ref{smoothness} it sandwiches the Hessian,
	\[
	\lambda_{\min}(X(\by))\,G_0\;\preceq\;\nabla^2f(\by)\;\preceq\;\tfrac12 G_0,
	\]
	the lower bound being attained at $\by=0$.
\end{remark}

\begin{lemma}[Strong convexity on sublevel sets]\label{strongconvexity}
	Let $\by_0\in\R^m$ and let $\mathscr S:=\{\by\in\R^m:f(\by)\le f(\by_0)\}$. Assume that
	$\bb\in\Int\mathscr M$, and define
	\begin{equation}\label{eq:delta}
	\delta:=\mathrm{dist}(\bb,\partial\mathscr M)
	=\min_{\|\bv\|_2=1}\bigl(\lambda_{\max}(A(\bv))-\bb^T\bv\bigr)>0 .
	\end{equation}
	Then $f$ is $\alpha$-strongly convex on $\mathscr S$, with
	\[
	\alpha:=\frac{\lambda_{\min}(G_0)}{n}
	\exp\!\left(-\frac{\sqrt{2\,\lambda_{\max}(G_0)}}{\delta}\,f(\by_0)\right)>0 .
	\]
\end{lemma}

\begin{proof}
	Let $h$ denote the support function of $\mathscr M$, so that $h(\bv)=\lambda_{\max}(A(\bv))$. Since
	$\bb\in\Int\mathscr M$, the support-function characterization of the interior gives
	$h(\bv)-\bb^T\bv>0$ for every $\bv\neq0$. By homogeneity, continuity and compactness of the unit
	sphere, $\delta=\min_{\|\bv\|_2=1}(h(\bv)-\bb^T\bv)>0$. This quantity is the radius of the largest
	Euclidean ball centered at $\bb$ contained in $\mathscr M$, hence equals
	$\mathrm{dist}(\bb,\partial\mathscr M)$.
	
	By Lemma~\ref{curvature}, for every $\by\in\R^m$,
	\[
	\nabla^2f(\by)\succeq\lambda_{\min}(X(\by))\,G_0\succeq\lambda_{\min}(X(\by))\,\lambda_{\min}(G_0)I_m ,
	\]
	where $G_0\succ0$ under Assumption~\ref{injective}. It remains to lower-bound
	$\lambda_{\min}(X(\by))$ uniformly on $\mathscr S$.
	
	Let $\mu_1\ge\cdots\ge\mu_n$ be the eigenvalues of $A(\by)$, and set $\gamma(\by):=\mu_1-\mu_n$.
	Since the eigenvalues of $X(\by)=\exp A(\by)/\Tr\exp A(\by)$ are $e^{\mu_k}/\sum_je^{\mu_j}$ and
	$\sum_je^{\mu_j}\le n\,e^{\mu_1}$,
	\[
	\lambda_{\min}(X(\by))\ge\frac1n e^{-\gamma(\by)} .
	\]
	The spread $\gamma(\by)$ is invariant under adding a scalar multiple of $I_n$ to $A(\by)$.
	Therefore, with $A^0(\by):=A(\by)-\tfrac{\Tr A(\by)}{n}I_n$, one obtains
	$\gamma(\by)^2\le2\Tr(A^0(\by)^2)=2\,\by^TG_0\by$, and hence
	\[
	\gamma(\by)\le\sqrt{2\,\lambda_{\max}(G_0)}\,\|\by\|_2 .
	\]
	On the other hand, $f(\by)=\log\Tr\exp A(\by)-\bb^T\by\ge h(\by)-\bb^T\by$. For $\by\neq0$, writing
	$\by=\|\by\|_2\bv$ with $\|\bv\|_2=1$, homogeneity gives
	$h(\by)-\bb^T\by=\|\by\|_2(h(\bv)-\bb^T\bv)\ge\delta\|\by\|_2$; the same inequality is trivial at
	$\by=0$. Thus $f(\by)\ge\delta\|\by\|_2$ for all $\by$, and in particular $f(\by_0)\ge0$. For
	$\by\in\mathscr S$ this yields $\|\by\|_2\le f(\by_0)/\delta$, so that
	\[
	\lambda_{\min}(X(\by))\ge\frac1n\exp\!\left(-\frac{\sqrt{2\,\lambda_{\max}(G_0)}}{\delta}\,f(\by_0)\right)
	\qquad(\by\in\mathscr S).
	\]
	Substituting into the Hessian lower bound gives $\nabla^2f(\by)\succeq\alpha I_m$ for every
	$\by\in\mathscr S$, with $\alpha$ as stated; hence $f$ is $\alpha$-strongly convex on $\mathscr S$.
\end{proof}
}

\section{Pre-conditioning}

{
For the conditioning $\kappa$ of $f$ to be as small as possible, Lemmas~\ref{smoothness}
and~\ref{strongconvexity} show that one should control both the upper and lower curvature, governed by the
Gram matrix $G_0$ of the traceless parts of the matrices $A_i$. In particular, small values of
$\lambda_{\max}(G_0)$ improve the smoothness bound, whereas large values of
$\lambda_{\min}(G_0)$ improve the strong convexity bound on sublevel sets.

This suggests a geometric preconditioning step: replace the original matrices $A_i$ by
a new set $\widehat A_i$, obtained by centering and orthonormalizing the traceless
directions, so that the corresponding centered Gram matrix is well conditioned, ideally $\widehat G_0=I_m$,
and both curvature estimates become optimally normalized. 
}

{Let us denote $\delta_{ij}:=1$ if $i=j$ and $0$ otherwise.}
	
	\begin{algorithm}[H]
		\caption{Pre-conditioning}\label{algo}
		\begin{algorithmic}[1]
			\REQUIRE \(A_1,\dots,A_m \in \S^n\)\\
			\ENSURE \(\hat A_1,\dots,\hat A_m \in \S^n\) with  \(\Tr\hat A_i=0\) and \( \Tr(\hat A_i\hat A_j)= \delta_{ij}\)
			\STATE {\bfseries (Center to traceless)}  
			\[
			{a_i:=\tfrac1n\Tr A_i,} \quad A'_i \;=\; A_i \;-a_i\,I_n,
			\quad i=1,\dots,m.
			\]
			\STATE {\bfseries (Compute Gram matrix)}  
			\[
			G'_{ij} \;=\;\Tr(A'_i\,A'_j),
			\quad i,j=1,\dots,m.
			\]
			\STATE {\bfseries (Whitening)}  
			Compute the symmetric inverse‐square‐root \(G'^{-1/2}\) via eigendecomposition:
			\[
			G' = U\,D\,U^T,
			\quad
			{W}:=G'^{-1/2} = U\,D^{-1/2}\,U^T.
			\]
			\STATE {\bfseries (Form Orthonormal Basis)}  
			\[
			\hat A_i
			=\sum_{j=1}^m {W_{ij}}\;A'_j,
			\quad i=1,\dots,m.
			\]
		\end{algorithmic}
	\end{algorithm}
	
	\begin{theorem}[Correctness]
		The output $\hat A_i$ of Algorithm \ref{algo} satisfies:
		\begin{enumerate}
			\item \(\displaystyle\Tr\hat A_i=0\) for all \(i\).
			\item  \( \Tr(\hat A_i\hat A_j)= \delta_{ij}\).
		\end{enumerate}
		Hence the \(\hat A_i\) are traceless and orthonormal.
	\end{theorem}
	
	\begin{proof}
		Each \(A'_i\) is by construction
		\[
		\Tr A'_i
		=\Tr A_i - {a_i} \Tr I_n
		=\Tr A_i - \Tr A_i
		=0.
		\]
		Since \(\hat A_i\) is a linear combination of the \(A'_j\), it too is traceless:
		\(\Tr \hat A_i=\sum_j {W_{ij}}\,\Tr A'_j=0.\)
		
		\medskip\noindent
		\textbf{2. Orthonormality.}  Define the centered Gram matrix \(G'\).  Then
		\[
		\Tr(\hat A_i\,\hat A_j)
		=\sum_{p,q}{W_{ip}}{W_{jq}}\,
		\Tr(A'_p\,A'_q)
		=\bigl[{W}\,G'\,{W}\bigr]_{ij}
		=\delta_{ij}.
		\]
		This shows \((\hat A_i)\) are orthonormal in the Frobenius inner product.
	\end{proof}

\begin{figure}[h]
	\begin{center}
		\includegraphics[width=0.9\textwidth]{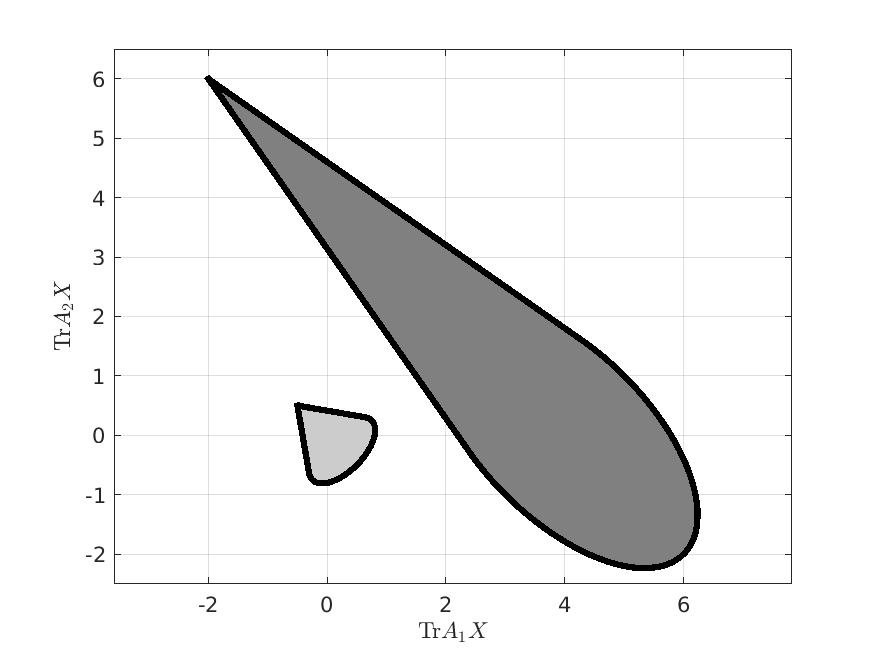}
		\caption{Moment bodies before (dark gray) and after (light gray) pre-conditioning Algorithm \ref{algo}.\label{fig:momentbodycond}}
	\end{center}
\end{figure}
\begin{example}
Let $n=3$, $m=2$ and
\[
A_1 = 
\left(\begin{array}{rrr}
6 & 1 & 0 \\ 1 & 2 & 0 \\ 0 & 0 & -2
\end{array}\right), \quad
A_2 = 
\left(\begin{array}{rrr}
-2 & 1 & 0 \\ 1 & 2 & 0 \\ 0 & 0 & 6
\end{array}\right)
\]
whose Gram matrix  has eigenvalues $28$ and $64$. Step 1 of Algorithm \ref{algo} yields the traceless matrices
\[
A'_1 = 
\left(\begin{array}{rrr}
4 & 1 & 0 \\ 1 & 0 & 0 \\ 0 & 0 & -4
\end{array}\right), \quad
A'_2 = 
\left(\begin{array}{rrr}
-4 & 1 & 0 \\ 1 & 0 & 0 \\ 0 & 0 & 4
\end{array}\right)
\]
whose Gram matrix computed in step 2 has eigenvalues $4$ and $64$. Finally, step 4 yields the traceless orthonormal matrices
\[
\hat{A}_1 = 
\frac{1}{2}\left(\begin{array}{rrr}
1 & 1 & 0 \\ 1 & 0 & 0 \\ 0 & 0 & -1
\end{array}\right), \quad
\hat{A}_2 = 
\frac{1}{2}\left(\begin{array}{rrr}
-1 & 1 & 0 \\ 1 & 0 & 0 \\ 0 & 0 & 1
\end{array}\right)
\]
corresponding to Example \ref{ex:2d}, whose Gram matrix is identity. The corresponding moment bodies, before and after the application of Algorithm \ref{algo}, are represented in Figure \ref{fig:momentbodycond}. The same pre-conditioned moment body is also represented on Figures \ref{fig:momentbody2d} and \ref{fig:momentbodygrid2d}.
\end{example}

We remark that the idea of using the Cholesky factor of the Gram matrix of the linear map was already exploited in the context of projection methods for semidefinite optimization, see \cite{mprw09} and \cite[section 3.2]{hm12}.

{The preconditioning changes coordinates in the moment body. Let $W:=G'^{-1/2}$ and $\ba:=(a_i)_i=(\tfrac1n\Tr A_i)_i$.   If \[ \widehat A_i:=\sum_{j=1}^m W_{ij}A'_j, \] then, for every $X\in\S^n_1$, \[ \widehat{\mathcal A}(X) = W(\mathcal A(X)-\ba). \] Consequently the original membership problem \[ \bb\in\mathscr M \] is equivalent to the preconditioned membership problem \[ \widehat\bb\in\widehat{\mathscr M}, \qquad \widehat\bb:=W(\bb-\ba), \qquad \widehat{\mathscr M}:=W(\mathscr M-\ba). \] Thus the underlying feasibility question is unchanged, but
	 both the moment coordinates and the target vector must be transformed. }
	 
{
The corresponding dual objective is \[ \widehat f(\widehat\by) := \log\Tr\exp\widehat A(\widehat\by) - \widehat\bb^T\widehat\by, \qquad \widehat A(\widehat\by):=\sum_{i=1}^m\widehat y_i\widehat A_i. \] Since \(W\) is symmetric, \[ \widehat A(\widehat\by) = A'(W\widehat\by), \qquad \widehat\bb=W(\bb-\ba), \] and therefore \[ \widehat f(\widehat\by) = f(W\widehat\by). \] Thus minimizers are related by \[ \by^\star=W\widehat\by^\star. \]
}
\section{Refined geometric analysis}

In this section we assume that the data are normalized: after the application of Algorithm \ref{algo} the matrices $A_i$ are traceless and orthonormal, i.e.
\[
\trace A_i = 0, \quad \trace (A_i A_j) = \delta_{ij}, \quad i,j=1,\ldots,m,
\]
or equivalently $\mathcal A(I_n)=0$, $\mathcal A \circ \mathcal A^T = I_m$. {The map
$\bu\mapsto A(\bu)$ is then a linear isometry from $\R^m$ onto the subspace $V:=\mathrm{span}(A_1,\dots,A_m)$ of the traceless symmetric matrices $\S^n_0$}. The corresponding moment body $\mathscr M$ is normalized\footnote{A Traceless Orthonormal Moment Body can be called a TOMB, evoking a solid, well-defined shape, where all the traceless moments rest.}, and we now report some of its geometric properties.

The radius of $\mathscr M$ is defined as follows
\[
\mathrm{rad}\,\mathscr M:=\max_{X\in\S^n_1}\|\mathcal A(X)\|_2 = \max_{\|\bu\|_2=1} \|\mathcal A(\bu\bu^T)\|_2.
\]
{
\begin{lemma}\label{radius}
	It holds
	\[
	\mathrm{rad}\,\mathscr M\;\le\;\sqrt{\tfrac{n-1}{n}}\;<\;1,
	\]
	with equality if and only if $V$ contains a matrix of the form $\bu\bu^T-\tfrac1n I_n$ for some unit
	vector $\bu$ (in particular when the $A_i$ span all of $\S^n_0$).
\end{lemma}

\begin{proof}
	Since the $A_i$ are traceless, $\mathcal A(I_n)=0$ and $\mathcal A(X)=\mathcal A(X-\tfrac1n I_n)$;
	since they are orthonormal, $\mathcal A$ is the coordinate map of the orthogonal projection $P_V$
	onto $V$, so $\|\mathcal A(X)\|_2=\|P_V(X-\tfrac1n I_n)\|_F\le\|X-\tfrac1n I_n\|_F$. Moreover
	\[
	\|X-\tfrac1n I_n\|_F^2=\|X\|_F^2-\tfrac2n\trace X+\tfrac1n\le 1-\tfrac2n+\tfrac1n=\tfrac{n-1}{n},
	\]
	using $0 \preceq X \preceq I_n \Rightarrow X^2 \preceq X \Rightarrow \|X\|_F^2=\trace(X^2)\le\trace X=1$. This proves $\mathrm{rad}\,\mathscr M\le\sqrt{(n-1)/n}$.
	The first inequality is an equality iff $X-\tfrac1n I_n\in V$, the second iff $X$ has rank one,
	i.e.\ $X=\bu\bu^T$; both hold simultaneously iff $\bu\bu^T-\tfrac1n I_n\in V$ for some unit $\bu$.
\end{proof}

Recall from \eqref{eq:support} that the support function of $\mathscr M$ is
$h(\bv)=\lambda_{\max}(A(\bv))$. The width of $\mathscr M$ in a unit direction $\bv$ is
$w(\bv):=h(\bv)+h(-\bv)$; since $A(-\bv)=-A(\bv)$ gives $h(-\bv)=-\lambda_{\min}(A(\bv))$,
\[
w(\bv)=\lambda_{\max}(A(\bv))-\lambda_{\min}(A(\bv))
\]
is the spectral spread of $A(\bv)$. The thickness and diameter of $\mathscr M$ are defined as follows
\[
\mathrm{thick}\,\mathscr M:=\min_{\|\bv\|_2=1}w(\bv), \qquad \mathrm{diam}\,\mathscr M:=\max_{\|\bv\|_2=1}w(\bv).
\]

\begin{lemma}\label{diamthick}
	It holds
	\[
	\sqrt{\frac{n}{\lfloor n/2\rfloor\,\lceil n/2\rceil}} 	=\begin{cases}\dfrac{2}{\sqrt n}&n=2k,\\[4pt]\sqrt{\dfrac{2k+1}{k(k+1)}}&n=2k+1\end{cases}
	\;\le\;\mathrm{thick}\,\mathscr M\;\le\;\mathrm{diam}\,\mathscr M\;\le\;\sqrt2 .
	\]
	The upper bound is attained iff $V$ contains a rank-two matrix with eigenvalues $\pm\tfrac1{\sqrt2}$;
	the lower bound is attained iff $V$ contains a Frobenius-unit traceless matrix whose eigenvalues take
	two values with multiplicities $\lfloor n/2\rfloor$ and $\lceil n/2\rceil$. In particular, if
	$A_1,\dots,A_m$ span all of $\S^n_0$ (i.e.\ $m=\tfrac{n(n+1)}2-1$), both equalities hold.
\end{lemma}

\begin{proof}
	Since $\bv\mapsto A(\bv)$ is an isometry of $\R^m$ onto $V$, as $\bv$ runs over the unit sphere
	$A(\bv)$ runs over the Frobenius-unit sphere of $V$, and $w(\bv)=\mathrm{sp}(A(\bv))$ with
	$\mathrm{sp}(M):=\lambda_{\max}(M)-\lambda_{\min}(M)$. Hence
	\[
	\mathrm{thick}\,\mathscr M=\min_{M\in V,\,\|M\|_F=1}\mathrm{sp}(M),\qquad
	\mathrm{diam}\,\mathscr M=\max_{M\in V,\,\|M\|_F=1}\mathrm{sp}(M),
	\]
	and since $V\subseteq\S^n_0$ it suffices to bound $\mathrm{sp}(M)$ over all traceless
	Frobenius-unit matrices. Let $\mu_1\ge\cdots\ge\mu_n$ be the eigenvalues of such an $M$, so
	$\sum_i\mu_i=0$, $\sum_i\mu_i^2=1$, and put $R:=\mu_1-\mu_n$.
	
	{Upper bound.} $R^2=(\mu_1-\mu_n)^2\le2(\mu_1^2+\mu_n^2)\le2\sum_i\mu_i^2=2$, so
	$\mathrm{sp}(M)\le\sqrt2$. The first inequality is an equality iff $\mu_1=-\mu_n$, the second iff
	$\mu_i=0$ for $i\ne1,n$; together with $\sum_i\mu_i^2=1$ this means $M$ has rank two with eigenvalues
	$\pm\tfrac1{\sqrt2}$.
	
	{Lower bound.} We prove that every traceless $M$ with eigenvalues $\mu_1\ge\cdots\ge\mu_n$ and
	range $R:=\mu_1-\mu_n$ satisfies $\sum_i\mu_i^2\le\frac{\lfloor n/2\rfloor\lceil n/2\rceil}{n}R^2$;
	together with $\|M\|_F^2=\sum_i\mu_i^2=1$ this gives
	$\mathrm{sp}(M)=R\ge\sqrt{n/(\lfloor n/2\rfloor\lceil n/2\rceil)}$. Equivalently, we bound the maximum
	of the convex function $\Phi(\bmu):=\sum_i\mu_i^2$ over the compact convex set
	$S_R:=\{\bmu\in\R^n:\sum_i\mu_i=0,\ \max_i\mu_i-\min_i\mu_i\le R\}$, and show it is attained at a
	two-valued spectrum.
	
	Let $\bmu$ maximize $\Phi$ over $S_R$, and set $a:=\max_i\mu_i$, $b:=\min_i\mu_i$. First $a-b=R$:
	otherwise scaling $\bmu$ by $R/(a-b)>1$ keeps it in $S_R$ and strictly increases $\Phi$. Second, at
	most one coordinate lies in the open interval $(b,a)$: if $\mu_i,\mu_j\in(b,a)$, then for small
	$\varepsilon$ the replacement $\mu_i\mapsto\mu_i+\varepsilon$, $\mu_j\mapsto\mu_j-\varepsilon$ stays in
	$S_R$ (it preserves $\sum_k\mu_k$ and leaves $a,b$ attained) and changes $\Phi$ by
	$2\varepsilon(\mu_i-\mu_j)+2\varepsilon^2$, which is positive for one sign of $\varepsilon$,
	contradicting maximality. Third, suppose exactly one coordinate equals $w\in(b,a)$, with $k$
	coordinates at $a$ and $n-1-k$ at $b$. Then $b=a-R$ and $\sum_i\mu_i=0$ make both $b$ and
	$w=-(n-1)a+(n-1-k)R$ affine in $a$, so $\Phi=ka^2+(n-1-k)b^2+w^2$ is a quadratic in $a$ with positive
	leading coefficient $n(n-1)$, hence convex; its maximum over the admissible interval (where
	$b\le w\le a$) is at an endpoint, i.e.\ where $w\in\{a,b\}$ and the spectrum is two-valued.
	
	Thus the maximum is attained with the eigenvalues taking only the values $a$ and $b$, say $k$ equal to
	$a$ and $n-k$ to $b$; then $\sum_i\mu_i=0$ forces $a=\frac{n-k}{n}R$, $b=-\frac{k}{n}R$, and
	\[
	\Phi=k\,a^2+(n-k)\,b^2=\frac{k(n-k)}{n}\,R^2\;\le\;\frac{\lfloor n/2\rfloor\lceil n/2\rceil}{n}\,R^2,
	\]
	maximizing the integer $k(n-k)$. Equality forces the spectrum to be two-valued with multiplicities
	$\lfloor n/2\rfloor$ and $\lceil n/2\rceil$.
	
	The bounds $\mathrm{thick}\,\mathscr M\ge\sqrt{n/(\lfloor n/2\rfloor\lceil n/2\rceil)}$ and
	$\mathrm{diam}\,\mathscr M\le\sqrt2$ follow, each attained iff $V$ contains the corresponding
	extremal matrix; when $V=\S^n_0$ both extremal matrices lie in $V$, giving the stated equalities.
\end{proof}

When $V=\S^n_0$ and \(n>2\), the full-space extremal values satisfy
\[
\mathrm{diam}\,\mathscr M=\sqrt2
<
2\sqrt{\frac{n-1}{n}},
\]
so the farthest point from the origin and the opposite farthest point cannot both be simultaneously attained as antipodal rank-one extremal images. This is consistent with the fact that the full traceless moment body is not centrally symmetric.
}

{
Recall that if $\bb \in \Int \mathscr M$, then the margin $\delta$ was defined in \eqref{eq:delta} of Lemma \ref{strongconvexity} as the radius of the largest Euclidean ball centered at $\bb$ and
contained in $\mathscr M$. Since $B(\bb,\delta)\subseteq\mathscr M\subseteq
B(0,\sqrt{(n-1)/n})$ by Lemma~\ref{radius},
\[
\delta\;\le\;\sqrt{\tfrac{n-1}{n}}-\|\bb\|_2\;<\;1-\|\bb\|_2,
\]
the margin is thus strictly smaller than the distance $1-\|\bb\|_2$ of $\bb$ to the unit sphere.
Let \[ \by^* := \arg \min_{\by \in \R^m} f(\by) \] denote the  minimizer of $f$, which is unique from Lemma \ref{unique}.
The geometric properties of $\mathscr M$ allow us to bound the value of $f$ at $\by^*$ and the norm
of $\by^*$ in terms of the margin of $\bb$.

\begin{lemma}\label{bounds}
	If $\bb \in \Int\mathscr M$, it holds:
	\begin{enumerate}
		\item[(i)] $0\leq f(\by^*)=-\trace(X^*\log X^*)\le\log n$, with $f(\by^*)=\log n$ if and only if $\bb=0$;
		\item[(ii)] $\displaystyle\|\by^*\|_2\;\le\;\frac{f(\by^*)}{\delta}\;\le\;\frac{\log n}{\delta}$.
	\end{enumerate}
\end{lemma}

\begin{proof}
	(i) By strong duality (Theorem~\ref{primaldual}),
	$f(\by^*)=\max\{-\trace(X\log X):X\in\S^n_1,\ \mathcal A(X)=\bb\}=-\trace(X^*\log X^*)$, the entropy of $X^*$. As $X^*=\exp_1 A(\by^*)\succ0$ has full rank, this entropy is strictly
	positive. For the upper bound, $f(\by^*)\le f(0)=\log\trace\exp(0)=\log n$, with equality iff
	$\by^*=0$, i.e.\ iff $\nabla f(0)=\mathcal A(\tfrac1n I_n)-\bb=-\bb$ vanishes, i.e.\ iff $\bb=0$.
	
	(ii) Since $\log\trace\exp M\ge\lambda_{\max}(M)$, for $\by\ne0$ with $\bv:=\by/\|\by\|_2$,
	\[
	f(\by)\ge\lambda_{\max}(A(\by))-\bb^T\by=\|\by\|_2(\lambda_{\max}(A(\bv))-\bb^T\bv)
	\ge\delta\,\|\by\|_2,
	\]
	which is the left inequality of \eqref{eq:raybound} and holds trivially at $\by=0$. Applying it at
	$\by^*$ and using (i) gives $\delta\,\|\by^*\|_2\le f(\by^*)\le\log n$.
\end{proof}
}

{
\begin{theorem}\label{smoothrefined}
	Function $f$ is $\tfrac12$-smooth.
\end{theorem}
\begin{proof}
	Since the $A_i$ are traceless and orthonormal, $A_i^0=A_i$, so the centered Gram matrix \eqref{gram0}
	is $G_0=[\Tr(A_iA_j)]_{i,j}=I_m$. Lemma~\ref{smoothness} then gives
	$\nabla^2 f(\by)\preceq\tfrac12 G_0=\tfrac12 I_m$ for all $\by$; in particular $\lambda_{\max}(G_0)=1$
	and $f$ is $\tfrac12$-smooth.
\end{proof}
}

{\begin{theorem}\label{strongconvexityrefined}
		If $\bb \in \Int\mathscr M$, the function $f$ is $\alpha$-strongly convex on the sublevel set $\{\by:f(\by)\le\log n\}$, with
		\[
		\alpha=\frac1n\,e^{-\sqrt2\,\log n/\delta}=n^{-(1+\sqrt2/\delta)}.
		\]
	\end{theorem}
	
	\begin{proof}
		Since the $A_i$ are traceless and orthonormal, $G_0=I_m$, so $\lambda_{\min}(G_0)=\lambda_{\max}(G_0)=1$.
		Taking $\by_0=0$ in Lemma~\ref{strongconvexity}, where $f(0)=\log n$, gives the sublevel modulus
		$$\alpha=\frac{\lambda_{\min}(G_0)}{n}\exp(-\sqrt{2\lambda_{\max}(G_0)}\,f(0)/\delta)
		=\frac1n e^{-\sqrt2\log n/\delta}.$$
\end{proof} }

\section{Convergence analysis}\label{convergence}

In this section we ask what minimizing \(f\) reveals about the position of \(\bb\) relative to the moment
body \(\mathscr M\), and how fast. Exactly one of three cases holds:
\[
\bb\in\Int\mathscr M\ \text{(strict feasibility)},\qquad
\bb\in\partial\mathscr M\ \text{(weak feasibility)},\qquad
\bb\notin\mathscr M\ \text{(infeasibility)} .
\]
Throughout we assume \(n\ge2\), the case \(n=1\) being trivial, and we work with the data normalized  by the preconditioning Algorithm~\ref{algo}. In particular,
by Theorem~\ref{smoothrefined}, \(f\) is \(\tfrac12\)-smooth.

The whole analysis rests on the gradient identity of Lemma~\ref{gh}, $\nabla f(\by)=\mathcal A(X(\by))-\bb$,
$X(\by)=\exp_1 A(\by)\in\S^n_1$,
in which every iterate furnishes both a point
$\mathcal A(X(\by))\in\Int\mathscr M$
and a direction \(\by\), the two ingredients of a certificate.

A single geometric quantity governs all three regimes: the distance
\[
\delta:=\mathrm{dist}(\bb,\partial\mathscr M)
\]
of \(\bb\) to the boundary of \(\mathscr M\). It is strictly positive precisely when
\(\bb\notin\partial\mathscr M\) and vanishes exactly on the boundary. In the strictly feasible case,
\(\delta\) is the radius of the largest Euclidean ball centered at \(\bb\) contained in
\(\mathscr M\), as in \eqref{eq:delta}. In the infeasible case it is the distance to the nearest point
of \(\mathscr M\), since for a closed convex set the metric projection of an exterior point lies on the
boundary, so that $\mathrm{dist}(\bb,\mathscr M)=\mathrm{dist}(\bb,\partial\mathscr M)$.

We shall use the following standard consequences of smoothness and strong convexity.

\begin{lemma}[Standard smooth-convex estimates]\label{stdsmoothconvex}
	Let \(f:\R^m\to\R\) be convex and \(\lambda\)-smooth, i.e.
	\[
	f(\bz)\le f(\by)+\nabla f(\by)^T(\bz-\by)
	+\frac{\lambda}{2}\|\bz-\by\|_2^2
	\qquad\forall\,\by,\bz\in\R^m.
	\]
	Let
	\[
	\by_{k+1}=\by_k-\frac1\lambda\nabla f(\by_k).
	\]
	Then the following estimates hold.
	
	\begin{enumerate}
		\item[(i)] Descent:
		\[
		f(\by_{k+1})
		\le
		f(\by_k)-\frac1{2\lambda}\|\nabla f(\by_k)\|_2^2 .
		\]
		
		\item[(ii)] Monotonicity of the gradient norm:
		\[
		\|\nabla f(\by_{k+1})\|_2
		\le
		\|\nabla f(\by_k)\|_2
		\qquad\forall k.
		\]
		
		\item[(iii)] Sublinear comparison estimate: for every \(\bz\in\R^m\) and every \(k\ge1\),
		\[
		f(\by_k)-f(\bz)
		\le
		\frac{\lambda\|\by_0-\bz\|_2^2}{2k}.
		\]
		
		\item[(iv)] Let $\mathscr C \subset \R^m$ be convex, let \(\by^\star \in \mathscr C\) be a minimizer of $f$, and suppose that \( f\) is
		\(\alpha\)-strongly convex on $\mathscr C$. Then
		for every $\by \in \mathscr C$,
		\[
		\|\nabla f(\by)\|_2^2
		\ge
		2\alpha(f(\by)-f(\by^\star)).
		\]
		
		\item[(v)] If \( f\) has a minimizer \(\by^\star\), then for every \(\by\in\R^m\),
		\[
		\|\nabla f(\by)\|_2^2
		\le
		2\lambda(f(\by)-f(\by^\star)).
		\]
	\end{enumerate}
\end{lemma}

\begin{proof}
	The descent estimate (i) follows by applying the quadratic upper bound with
	\(\bz=\by-\lambda^{-1}\nabla f(\by)\). This is the standard descent estimate for smooth functions \cite[Sec.~2.2.4]{Nesterov18}.
	
	For (ii), we first recall a standard consequence of smooth convexity. For every \(\bu,\bv\in\R^m\),
	\begin{equation}\label{eq:grad-firm}
		\frac1\lambda
		\|\nabla f(\bu)-\nabla f(\bv)\|_2^2
		\le
		\langle\nabla f(\bu)-\nabla f(\bv),\bu-\bv\rangle
	\end{equation}
	see \cite[Thm.~2.1.5, eq.~(2.1.11)]{Nesterov18}.	
	Now apply \eqref{eq:grad-firm} with \(\bu=\by_{k+1}\) and \(\bv=\by_k\). Since $\by_{k+1}-\by_k=-\frac1\lambda\nabla f(\by_k)$, we obtain
	\[
	-\left\langle
	\nabla f(\by_{k+1})-\nabla f(\by_k),
	\nabla f(\by_k)
	\right\rangle
	\ge
	\|\nabla f(\by_{k+1})-\nabla f(\by_k)\|_2^2.
	\]
	Expanding the right-hand side gives
	\[
	-\langle\nabla f(\by_{k+1}),\nabla f(\by_k)\rangle
	+\|\nabla f(\by_k)\|_2^2
	\ge
	\|\nabla f(\by_{k+1})\|_2^2
	-2\langle\nabla f(\by_{k+1}),\nabla f(\by_k)\rangle
	+\|\nabla f(\by_k)\|_2^2.
	\]
	After cancellation, this becomes
	\[
	\|\nabla f(\by_{k+1})\|_2^2
	\le
	\langle\nabla f(\by_{k+1}),\nabla f(\by_k)\rangle.
	\]
	By Cauchy-Schwarz,
	\[
	\|\nabla f(\by_{k+1})\|_2^2
	\le
	\|\nabla f(\by_{k+1})\|_2\,\|\nabla f(\by_k)\|_2.
	\]
	If \(\nabla f(\by_{k+1})=0\), the claim is immediate. Otherwise we divide by
	\(\|\nabla f(\by_{k+1})\|_2\), and obtain
	\[
	\|\nabla f(\by_{k+1})\|_2
	\le
	\|\nabla f(\by_k)\|_2.
	\]
	
	For (iii), convexity gives
	\[
	f(\by_j)-f(\bz)
	\le
	\nabla f(\by_j)^T(\by_j-\bz).
	\]
	Combining this with the descent estimate (i), we obtain
	\[
	f(\by_{j+1})-f(\bz)
	\le
	\nabla f(\by_j)^T(\by_j-\bz)
	-\frac1{2\lambda}\|\nabla f(\by_j)\|_2^2.
	\]
	Using $\by_{j+1}=\by_j-\frac1\lambda\nabla f(\by_j)$, for every \(\bz\),
	\[
	f(\by_{j+1})-f(\bz)
	\le
	\frac{\lambda}{2}
	\left(
	\|\by_j-\bz\|_2^2-\|\by_{j+1}-\bz\|_2^2
	\right).
	\]
	Summing from \(j=0\) to \(k-1\) yields
	\[
	\sum_{j=0}^{k-1}
	(f(\by_{j+1})-f(\bz))
	\le
	\frac{\lambda}{2}\|\by_0-\bz\|_2^2.
	\]
	Since \(f(\by_k)\) is non-increasing by (i), we have
	\[
	k(f(\by_k)-f(\bz))
	\le
	\sum_{j=0}^{k-1}
	(f(\by_{j+1})-f(\bz)),
	\]
	which proves (iii). This is the standard \(O(1/k)\) gradient-method comparison estimate; compare
	with the analysis in \cite[Thm.~2.1.14 and Cor.~2.1.2]{Nesterov18}.
	
	For (iv), strong convexity on the segment \([\by,\by^\star]\) gives
	\[
	f(\by^\star)
	\ge
	f(\by)+\nabla f(\by)^T(\by^\star-\by)
	+\frac{\alpha}{2}\|\by-\by^\star\|_2^2.
	\]
	Hence
	\[
	f(\by)-f(\by^\star)
	\le
	\nabla f(\by)^T(\by-\by^\star)
	-\frac{\alpha}{2}\|\by-\by^\star\|_2^2
	\le
	\frac1{2\alpha}\|\nabla f(\by)\|_2^2,
	\]
	where the last inequality is the elementary bound
	$\ba^T\bu-\frac{\alpha}{2}\|\bu\|_2^2
	\le \frac1{2\alpha}\|\ba\|_2^2$.
	
	Finally, (v) follows from the descent estimate (i) and optimality of \(\by^\star\). Indeed,
	\[
	f(\by^\star)
	\le
	f\left(\by-\frac1\lambda\nabla f(\by)\right)
	\le
	f(\by)-\frac1{2\lambda}\|\nabla f(\by)\|_2^2.
	\]
	Rearranging gives
	\[
	\|\nabla f(\by)\|_2^2
	\le
	2\lambda(f(\by)-f(\by^\star)).
	\]
\end{proof}

Since \(f\) is \(\tfrac12\)-smooth, the gradient method with fixed step \(2\) is
\begin{equation}\label{eq:gradmethod}
	\by_{k+1}=\by_k-2\nabla f(\by_k).
\end{equation}
Applying Lemma~\ref{stdsmoothconvex}(i) gives the descent inequality
\begin{equation}\label{eq:descent}
	f(\by_{k+1})
	\le
	f(\by_k)-\|\nabla f(\by_k)\|_2^2.
\end{equation}

The three regimes are summarized as follows:
\[
\begin{array}{llll}
	\toprule
	& \by_k & \|\nabla f(\by_k)\|_2 & f(\by_k) \\
	\midrule
	\bb\in\Int\mathscr M
	& \to\by^\star
	& \to0\ \text{geometrically}
	& \to f^\star\in(0,\log n] \\
	\bb\in\partial\mathscr M
	& \|\by_k\|_2\to\infty
	& \to0,\ O(k^{-1/2})
	& \downarrow\inf f\in[0,\log n] \\
	\bb\notin\mathscr M
	& \|\by_k\|_2\to\infty
	& \ge\delta>0
	& \le f(\by_0)-k\delta^2\to-\infty \\
	\bottomrule
\end{array}
\]
The rigorous content is Theorems~\ref{strictlyfeasible}, \ref{weaklyfeasible} and
\ref{infeasiblecase}; the decision itself rests on the two certificates of the next subsection.

\subsection{Certificates of membership and distance}

The sign of \(f\) decides membership, and a single iterate brackets the distance to \(\mathscr M\) from both
sides.

\begin{lemma}\label{valuecert}
	\(\bb\in\mathscr M\) if and only if
	\[
	f(\by)\ge0\qquad\forall\,\by\in\R^m.
	\]
	Consequently, if \(f(\by_k)<0\) at some iterate, then \(\bb\notin\mathscr M\), and
	\(\bv:=\by_k\) separates \(\bb\) from \(\mathscr M\), namely
	\[
	\lambda_{\max}(A(\bv))<\bb^T\bv.
	\]
\end{lemma}

\begin{proof}
	If \(\bb\in\mathscr M\), then
	\[
	h(\by)=\lambda_{\max}(A(\by))
	=\max_{\bx\in\mathscr M}\bx^T\by
	\ge
	\bb^T\by
	\]
	by \eqref{eq:support}. Also,
	\[
	\log\trace\exp A(\by)
	\ge
	\lambda_{\max}(A(\by)).
	\]
	Therefore
	\[
	f(\by)
	\ge
	h(\by)-\bb^T\by
	\ge0.
	\]
	Conversely, if \(\bb\notin\mathscr M\), then by separation there exists \(\bv\in\R^m\) such that
	\[
	\bb^T\bv>h(\bv)=\lambda_{\max}(A(\bv)).
	\]
	Along the ray \(t\bv\),
	\[
	f(t\bv)
	\le
	t\lambda_{\max}(A(\bv))+\log n-t\bb^T\bv
	\to-\infty
	\qquad(t\to+\infty).
	\]
	Hence \(f\) takes negative values.
	
	Finally, if \(f(\by_k)<0\), then
	\[
	\lambda_{\max}(A(\by_k))-\bb^T\by_k
	\le
	f(\by_k)
	<0,
	\]
	which gives the separating inequality.
\end{proof}

\begin{lemma}\label{distbounds}
	For every \(\by\ne0\), with
	\[
	\bv:=\frac{\by}{\|\by\|_2}
	\qquad\text{and}\qquad
	X(\by)=\exp_1 A(\by),
	\]
	one has
	\[
	\bigl[\bb^T\bv-\lambda_{\max}(A(\bv))\bigr]_+
	\le
	\mathrm{dist}(\bb,\mathscr M)
	\le
	\|\mathcal A(X(\by))-\bb\|_2
	=
	\|\nabla f(\by)\|_2 .
	\]
\end{lemma}

\begin{proof}
	For the upper bound, \(\mathcal A(X(\by))\in\mathscr M\), and therefore
	\[
	\mathrm{dist}(\bb,\mathscr M)
	\le
	\|\mathcal A(X(\by))-\bb\|_2
	=
	\|\nabla f(\by)\|_2.
	\]
	For the lower bound, let \(\bx\in\mathscr M\). Since \(\|\bv\|_2=1\),
	\[
	\|\bb-\bx\|_2
	\ge
	\bv^T(\bb-\bx)
	=
	\bb^T\bv-\bx^T\bv
	\ge
	\bb^T\bv-h(\bv)
	=
	\bb^T\bv-\lambda_{\max}(A(\bv)).
	\]
	Minimizing over \(\bx\in\mathscr M\) gives
	\[
	\mathrm{dist}(\bb,\mathscr M)
	\ge
	\bb^T\bv-\lambda_{\max}(A(\bv)).
	\]
	Since the distance is nonnegative, the positive part may be taken.
\end{proof}

Both quantities are computable from a single iterate \(\by_k\): the upper bound
\(\|\nabla f(\by_k)\|_2\) certifies
\[
\mathrm{dist}(\bb,\mathscr M)\le\|\nabla f(\by_k)\|_2,
\]
while a strictly positive value of
\[
\bb^T\bv_k-\lambda_{\max}(A(\bv_k)),
\qquad
\bv_k:=\frac{\by_k}{\|\by_k\|_2},
\]
certifies \(\bb\notin\mathscr M\), together with a separating direction. The upper bound, read as
\[
\|\nabla f(\by)\|_2\ge\mathrm{dist}(\bb,\mathscr M)
\qquad\forall\,\by,
\]
is what drives the infeasible case below.

\subsection{Convergence of the gradient method}

We treat the three regimes in turn. In the strictly feasible case the iterates converge geometrically
to the unique minimizer. For any feasible \(\bb\), the residual still vanishes algebraically even
without using a positive margin; on the boundary \(\bb\in\partial\mathscr M\), where no minimizer
exists, this becomes the operative guarantee. In the infeasible case the value diverges to
\(-\infty\) at least linearly in \(k\), yielding a separating certificate after finitely many steps.

\begin{theorem}[strict feasibility]\label{strictlyfeasible}
	Let \(\bb\in\Int\mathscr M\), let
	\[
	\delta:=\mathrm{dist}(\bb,\partial\mathscr M)>0
	\]
	as in \eqref{eq:delta}, and let \((\by_k)\) be generated by \eqref{eq:gradmethod} from \(\by_0\).
	Given \(\epsilon>0\), the bound
	\[
	\|\nabla f(\by_k)\|_2\le\epsilon
	\]
	holds whenever
	\[
	k\ge
	n\,e^{\sqrt2\,f(\by_0)/\delta}
	\log\frac{\sqrt{f(\by_0)}}{\epsilon}.
	\]
	In particular, for \(\by_0=0\), where \(f(0)=\log n\),
	\[
	k\ge
	n^{1+\sqrt2/\delta}
	\log\frac{\sqrt{\log n}}{\epsilon}.
	\]
\end{theorem}

\begin{proof}
	By Theorem~\ref{smoothrefined}, \(f\) is \(\lambda\)-smooth with \(\lambda=\tfrac12\). Since the data
	are normalized, \(G_0=I_m\). Applying Lemma~\ref{strongconvexity} to the sublevel set
	\[
	\mathscr S_0:=\{\by:f(\by)\le f(\by_0)\}
	\]
	shows that \(f\) is \(\alpha\)-strongly convex on \(\mathscr S_0\), with
	\[
	\alpha=
	\frac1n
	\exp\left(-\frac{\sqrt2\,f(\by_0)}{\delta}\right).
	\]
	In particular, since \(n\ge2\),
	\[
	\alpha\le\frac1n\le\lambda.
	\]
	
	By Lemma~\ref{stdsmoothconvex}(i), or equivalently by \eqref{eq:descent}, the values decrease.
	Hence all iterates remain in the convex set \(\mathscr S_0\). The minimizer \(\by^\star\) also belongs
	to \(\mathscr S_0\), since \(f(\by^\star)\le f(\by_0)\).
	
	For every \(\by\in\mathscr S_0\), the segment \([\by,\by^\star]\) is contained in \(\mathscr S_0\).
	Therefore Lemma~\ref{stdsmoothconvex}(iv) gives
	\[
	\|\nabla f(\by)\|_2^2
	\ge
	2\alpha\bigl(f(\by)-f(\by^\star)\bigr).
	\]
	Combining this with the descent estimate yields
	\[
	f(\by_{k+1})-f(\by^\star)
	\le
	\left(1-\frac{\alpha}{\lambda}\right)
	\bigl(f(\by_k)-f(\by^\star)\bigr).
	\]
	Thus, by induction and \(1-t\le e^{-t}\),
	\[
	f(\by_k)-f(\by^\star)
	\le
	e^{-(\alpha/\lambda)k}
	\bigl(f(\by_0)-f(\by^\star)\bigr).
	\]
	Using \(f(\by^\star)\ge0\), this gives
	\[
	f(\by_k)-f(\by^\star)
	\le
	e^{-(\alpha/\lambda)k}f(\by_0).
	\]
	
	Finally, Lemma~\ref{stdsmoothconvex}(v) gives
	\[
	\|\nabla f(\by_k)\|_2
	\le
	\sqrt{2\lambda f(\by_0)}
	e^{-(\alpha/2\lambda)k}.
	\]
	Requiring the right-hand side to be at most \(\epsilon\) gives
	\[
	k
	\ge
	\frac{2\lambda}{\alpha}
	\log\frac{\sqrt{2\lambda f(\by_0)}}{\epsilon}.
	\]
	Since \(\lambda=\tfrac12\),
	\[
	k
	\ge
	\frac1\alpha
	\log\frac{\sqrt{f(\by_0)}}{\epsilon}
	=
	n\,e^{\sqrt2 f(\by_0)/\delta}
	\log\frac{\sqrt{f(\by_0)}}{\epsilon}.
	\]
	For \(\by_0=0\), \(f(0)=\log n\) and
	\[
	\frac1\alpha=n^{1+\sqrt2/\delta}.
	\]
	
	The same contraction also gives \(\by_k\to\by^\star\). Indeed, strong convexity on
	\(\mathscr S_0\) yields
	\[
	\frac{\alpha}{2}\|\by_k-\by^\star\|_2^2
	\le
	f(\by_k)-f(\by^\star).
	\]
\end{proof}

\begin{theorem}[feasible case, no margin]\label{weaklyfeasible}
	Let \(\bb\in\mathscr M\), and let \((\by_k)\) be generated by \eqref{eq:gradmethod} from \(\by_0\).
	Then
	\[
	\sum_{j=0}^{\infty}\|\nabla f(\by_j)\|_2^2
	\le
	f(\by_0),
	\]
	the residual norm is non-increasing, and
	\[
	\|\nabla f(\by_k)\|_2
	\le
	\sqrt{\frac{f(\by_0)}{k+1}},
	\qquad
	f(\by_k)\downarrow \inf_{\by}f(\by)\ge0.
	\]
	Equivalently,
	\[
	\mathcal A(X(\by_k))-\bb\to0.
	\]
	When \(\bb\in\partial\mathscr M\), the function \(f\) has no minimizer and
	\[
	\|\by_k\|_2\to\infty.
	\]
\end{theorem}

\begin{proof}
	Since \(\bb\in\mathscr M\), Lemma~\ref{valuecert} gives \(f\ge0\) on \(\R^m\). Summing
	\eqref{eq:descent} from \(j=0\) to \(k\) gives
	\[
	\sum_{j=0}^{k}\|\nabla f(\by_j)\|_2^2
	\le
	f(\by_0)-f(\by_{k+1})
	\le
	f(\by_0).
	\]
	Letting \(k\to\infty\) gives
	\[
	\sum_{j=0}^{\infty}\|\nabla f(\by_j)\|_2^2
	\le
	f(\by_0),
	\]
	and hence
	\[
	\|\nabla f(\by_k)\|_2\to0.
	\]
	
	By Lemma~\ref{stdsmoothconvex}(ii), the residual norm is non-increasing. Therefore
	\[
	\|\nabla f(\by_k)\|_2^2
	=
	\min_{0\le j\le k}
	\|\nabla f(\by_j)\|_2^2
	\le
	\frac1{k+1}
	\sum_{j=0}^{k}\|\nabla f(\by_j)\|_2^2
	\le
	\frac{f(\by_0)}{k+1}.
	\]
	This proves the last-iterate \(O(k^{-1/2})\) residual bound.
	
	For value convergence, Lemma~\ref{stdsmoothconvex}(iii) gives, for every \(\bz\in\R^m\) and every
	\(k\ge1\),
	\[
	f(\by_k)-f(\bz)
	\le
	\frac{\lambda\|\by_0-\bz\|_2^2}{2k}.
	\]
	Letting \(k\to\infty\) gives
	\[
	\limsup_{k\to\infty} f(\by_k)\le f(\bz)
	\qquad\forall\,\bz\in\R^m.
	\]
	Hence
	\[
	\limsup_{k\to\infty} f(\by_k)\le \inf_{\by}f(\by).
	\]
	Since \(f(\by_k)\) is decreasing and bounded below by \(0\), it follows that
	\[
	f(\by_k)\downarrow\inf_{\by}f(\by).
	\]
	
	Finally, assume that \(\bb\in\partial\mathscr M\). Then \(\nabla f\) never vanishes. Indeed, if
	\(\nabla f(\by)=0\), then $\mathcal A(X(\by))=\bb$,
	and by Lemma~\ref{diffeo} this point belongs to \(\Int\mathscr M\), contradicting
	\(\bb\in\partial\mathscr M\). Since \(\|\nabla f(\by_k)\|_2\to0\), any bounded subsequence of
	\((\by_k)\) would have a convergent subsequence accumulating at a zero of \(\nabla f\). This is
	impossible. Therefore
	\[
	\|\by_k\|_2\to\infty.
	\]
\end{proof}

\begin{theorem}[infeasibility]\label{infeasiblecase}
	Let \(\bb\notin\mathscr M\),
	and let \((\by_k)\) be generated by \eqref{eq:gradmethod} from \(\by_0\). Then
	\[
	f(\by_k)
	\le
	f(\by_0)-k\delta^2.
	\]
	Consequently, \(f(\by_k)<0\) for every integer
	\[
	k>\frac{f(\by_0)}{\delta^2}.
	\]
	In particular, if \(\by_0=0\), then \(f(\by_k)<0\) for every integer
	\[
	k>\frac{\log n}{\delta^2}.
	\]
	At any iterate with \(f(\by_k)<0\), the vector \(\bv:=\by_k\) separates \(\bb\) from
	\(\mathscr M\):
	\[
	\lambda_{\max}(A(\bv))<\bb^T\bv.
	\]
	Thus infeasibility is certified. Moreover,
	\[
	f(\by_k)\to-\infty,
	\qquad
	\|\by_k\|_2\to\infty.
	\]
\end{theorem}

\begin{proof}
	Since \(\mathscr M\) is closed and convex and \(\bb\notin\mathscr M\), the metric projection of
	\(\bb\) onto \(\mathscr M\) lies on \(\partial\mathscr M\). Hence
	\[
	\mathrm{dist}(\bb,\mathscr M)
	=
	\mathrm{dist}(\bb,\partial\mathscr M)
	=
	\delta.
	\]
	By Lemma~\ref{distbounds},
	\[
	\|\nabla f(\by)\|_2
	=
	\|\mathcal A(X(\by))-\bb\|_2
	\ge
	\delta
	\qquad\forall\,\by.
	\]
	Using \eqref{eq:descent}, we get
	\[
	f(\by_{k+1})
	\le
	f(\by_k)-\|\nabla f(\by_k)\|_2^2
	\le
	f(\by_k)-\delta^2.
	\]
	Iterating gives
	\[
	f(\by_k)\le f(\by_0)-k\delta^2.
	\]
	This is negative whenever
	\[
	k>\frac{f(\by_0)}{\delta^2}.
	\]
	The separating certificate then follows from Lemma~\ref{valuecert}. Finally, since
	\(f(\by_k)\to-\infty\), the sequence \((\by_k)\) cannot be bounded, because \(f\) is finite and
	continuous on bounded sets. Hence
	\[
	\|\by_k\|_2\to\infty.
	\]
\end{proof}

The three rates make the trichotomy quantitative and complementary. A vanishing residual certifies
membership: geometrically in the strictly feasible case, and algebraically in the feasible case without
a margin. Conversely, when \(\bb\notin\mathscr M\), the value becomes negative after every integer
\[
k>\frac{f(\by_0)}{\delta^2},
\]
and the corresponding iterate provides an explicit separating hyperplane.

\subsection{Sufficient conditions from the data}

For the normalized problem, the norm of $\bb$ alone settles two of the three cases, through the
inscribed and circumscribed central balls.

\begin{proposition}\label{ballconditions}
	Assume the data are normalized ($\trace A_i=0$, $\trace(A_iA_j)=\delta_{ij}$). Then
	\[
	B\Bigl(0,\tfrac{1}{\sqrt{n(n-1)}}\Bigr)\;\subseteq\;\mathscr M\;\subseteq\;
	B\Bigl(0,\sqrt{\tfrac{n-1}{n}}\Bigr),
	\]
	so that
	\begin{enumerate}
		\item[(i)] $\|\bb\|_2<\tfrac{1}{\sqrt{n(n-1)}}\ \Rightarrow\ \bb\in\Int\mathscr M$;
		\item[(ii)] $\|\bb\|_2>\sqrt{\tfrac{n-1}{n}}\ \Rightarrow\ \bb\notin\mathscr M$.
	\end{enumerate}
\end{proposition}

\begin{proof}
	The outer inclusion is Lemma~\ref{radius}. For the inner one, fix a unit $\bv$ and let
	$t:=\lambda_{\max}(A(\bv))$ with eigenvalues $\mu_1,\dots,\mu_n$ of $A(\bv)$. Each $\mu_i\le t$, and
	$\sum_j\mu_j=0$ forces $\mu_i=-\sum_{j\ne i}\mu_j\ge-(n-1)t$; thus $\mu_i\in[-(n-1)t,t]$, so
	$(\mu_i-t)(\mu_i+(n-1)t)\le0$, i.e.\ $\mu_i^2\le-(n-2)t\,\mu_i+(n-1)t^2$. Summing and using
	$\sum_i\mu_i=0$,
	\[
	\|\bv\|_2^2=\|A(\bv)\|_F^2=\sum_i\mu_i^2\le n(n-1)\,t^2,
	\]
	hence $\lambda_{\max}(A(\bv))=t\ge\tfrac{1}{\sqrt{n(n-1)}}$. Therefore
	$h(\bv)\ge\tfrac{1}{\sqrt{n(n-1)}}$ for every unit $\bv$, which is the inner inclusion.
	
	(i) If $\|\bb\|_2<\tfrac{1}{\sqrt{n(n-1)}}$ then for all unit $\bv$,
	$\bb^T\bv\le\|\bb\|_2<\tfrac{1}{\sqrt{n(n-1)}}\le\lambda_{\max}(A(\bv))$, so $g(\bv)>0$ and
	$\bb\in\Int\mathscr M$ by Lemma~\ref{coercive}. (ii) If $\|\bb\|_2>\sqrt{(n-1)/n}\ge\mathrm{rad}\,\mathscr M$
	then 
	$\bb\notin\mathscr M$.
\end{proof}

\section{Block-separable problems}

A natural structured instance arises when the constraint matrices are block diagonal, i.e.\ direct
sums $A_i=\bigoplus_{j=1}^p A_{ij}$ with blocks $A_{ij}\in\S^{n_j}$, $i=1,\dots,m$. Let
$N:=\sum_{j=1}^p n_j$ denote the total size. Because each $A_i$ is block diagonal, $\Tr(A_iX)$ depends
only on the block-diagonal part of $X$, and replacing $X$ by its block-diagonal part does
not decrease the von Neumann entropy; hence the maximum-entropy matrix consistent with the data is
itself block diagonal. The maximum-entropy primal of Theorem~\ref{primaldual} may therefore be
restricted, without loss, to block-diagonal densities $X=\bigoplus_{j=1}^p X_j$:
\[
\begin{aligned}
	\max_{X_1,\dots,X_p}\;& -\sum_{j=1}^p\Tr(X_j\log X_j)\\
	\text{s.t.}\quad
	& \sum_{j=1}^p\Tr(A_{ij}X_j)=b_i,\quad i=1,\dots,m,\\
	& X_j\succeq0,\quad \sum_{j=1}^p\Tr(X_j)=1 .
\end{aligned}
\]
Assume $\bb$ lies in the interior of the corresponding moment body $\mathscr M=\mathcal A(\S^N_1)$, so
that the Slater point $X=\tfrac1N I_N$ certifies strict feasibility. The Lagrangian separates over
$j$, and exactly as in the proof of Theorem~\ref{primaldual} strong duality gives the unconstrained
dual
\[
\min_{\by\in\R^m}\;
\log\sum_{j=1}^p \Tr\exp\Bigl(\sum_{i=1}^m y_i\,A_{ij}\Bigr)\;-\;\bb^\top\by ,
\]
the only change from the single-block case being that the exponent $A(\by)=\bigoplus_j\sum_i y_iA_{ij}$
is block diagonal, so that $\Tr\exp A(\by)=\sum_{j}\Tr\exp(\sum_i y_iA_{ij})$. The unique primal optimum
is the jointly normalized block Gibbs state
\[
X_j^\star=\frac{\exp\bigl(\sum_i y_i^\star A_{ij}\bigr)}{Z(\by^\star)},\qquad
Z(\by)=\sum_{j=1}^p\Tr\exp\Bigl(\sum_i y_i A_{ij}\Bigr).
\]

The benefit is computational. Evaluating the dual objective and its gradient now requires only the $p$
small matrix exponentials $\exp(\sum_i y_iA_{ij})$ of sizes $n_j$, at total cost $\sum_j O(n_j^3)$
rather than the $O(N^3)$ of a single $N\times N$ exponential, while the curvature estimates and the
preconditioning strategy apply unchanged: the centered Gram matrix
$G_0=[\Tr(A_i^0A_j^0)]_{ij}$ assembles additively over blocks, and the dual remains $\tfrac12$-smooth
after normalization.

\section{Numerical experiments}

We constructed a basic Matlab implementation\footnote{Available for download at \hyperlink{https://homepages.laas.fr/henrion/software/maxentmom/maxentmom.m}{homepages.laas.fr/henrion/software/maxentmom/maxentmom.m}} of L-BFGS that takes as input matrix $A$ of size $m$-by-$n^2$ and a vector $\bb$ of size $m$, and returns a vector $\by$ of size $m$ minimizing $f$:
\begin{verbatim}
y = maxentmom(A,b);
\end{verbatim}
 The algorithm calls the following function which evaluates $f$ and its gradient:
\begin{verbatim}
function [val, grad] = logpart(A,b,y)
% A : matrix of size m by n^2
% b, y : vectors of size m
n = sqrt(size(A,2));
[V,D] = eig(reshape(A'*y,n,n));
X = V * diag(exp(diag(D))) * V';
t = trace(X); 
val = log(t) - y'*b; % f(y)
grad = (A*X(:))/t - b; % grad f(y)
end
\end{verbatim}
Alternatively we can use HANSO \cite{hanso} which is a Matlab implementation of L-BFGS also aimed at non-smooth non-convex problems.

Convergence of iterate $\by_k$ occurs when the norm of the residual $\mathcal{A}(\exp_1 A(\by_k))-\bb$ (i.e. the gradient of $f$ at $\by_k$) is smaller than some a priori given expected accuracy, typically $10^{-8}$.
This is a relative accuracy when the data is normalized via Algorithm \ref{algo}, since the norm of $\bb$ is less than one by Lemma \ref{radius}.

{	All experiments were carried out in Matlab \texttt{R2025a} on a standard laptop
running Ubuntu~24.04, equipped with an Intel Core i7-1165G7 processor
(4 physical cores, 8 threads) at 2.8~GHz and 16~GB of memory. For comparison we used the two state-of-the-art semidefinite solvers SDPNAL+
version~1.0 \cite{sdpnal15,sdpnal20} and  MOSEK 
version~11.0.25 \cite{mosek}, the latter called through its native
\texttt{mosekopt} interface. Reported timings are wall-clock times.
}

%\subsection{Toy problem}
%\begin{figure}[h]
%	\begin{center}
%		\includegraphics[width=0.8\textwidth]{iterates2d.jpg}
%		\caption{Typical iterates (black dots) starting from the origin and reaching various target vectors near the boundary (dark gray) of the moment body (light gray).\label{fig:iterates2d}}
%	\end{center}
%\end{figure}
%
%Let us illustrate the behavior of {\tt maxentmom} on our toy planar moment body of Example \ref{ex:2d}. On Figure \ref{fig:iterates2d} are represented 10 typical trajectories $\mathcal A(\exp_1 A(\by_k))$ for 10 different target vectors $\bb$ chosen close to the boundary of the moment body, with the same initial condition $\by_0 = 0$. Iterates are represented by black dots, and typically 7 iterations suffice to reach the target vector at accuracy $10^{-8}$. 

\subsection{Entropy landscape and Gibbs spectrum}
Staying with the planar moment body of Example~\ref{ex:2d}, we visualize the two
objects underlying the duality of Theorem~\ref{primaldual}: the optimal value of $f$
as a function of the target $\bb$, and the Gibbs state $X(\by)$ along a ray.

\begin{figure}[h]
	\begin{center}
		\includegraphics[width=0.72\textwidth]{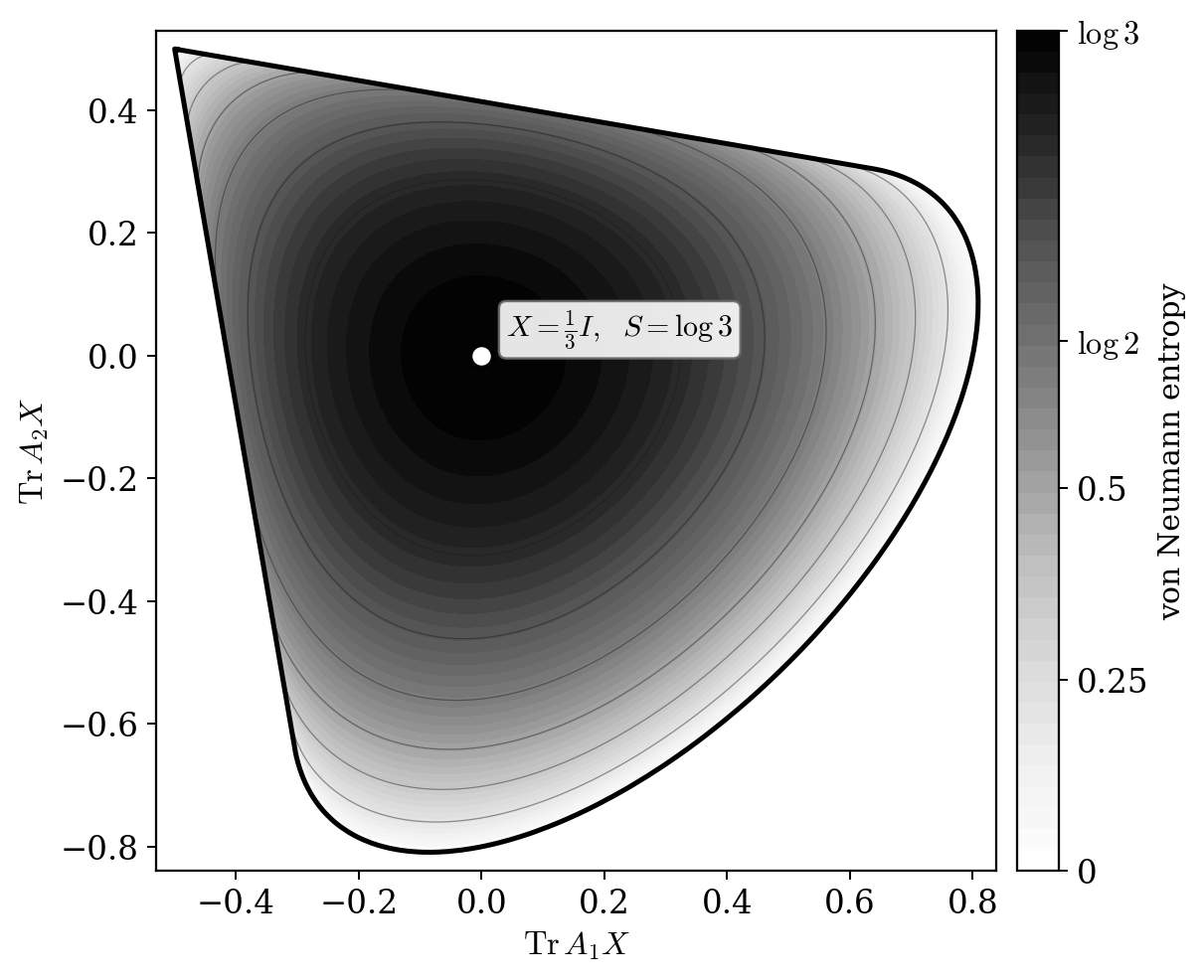}
		\caption{Maximal von Neumann entropy $S(\bb)=\min_{\by} f(\by)
			=-\Tr(X^\star\log X^\star)$ over the moment body $\mathscr M$ of
			Example~\ref{ex:2d} (darker gray is higher entropy). It peaks at $\log 3$ at the
			center $\bb=\mathbf 0$ (maximally mixed $X^\star=\tfrac13 I$, white marker),
			decays to $0$ on the curved boundary (rank-one $X^\star$), and stays elevated up
			to $\log 2$ on the two flat edges (rank-two $X^\star$). Thin lines are entropy
			level sets.\label{fig:entropy_landscape}}
	\end{center}
\end{figure}

By Theorem~\ref{primaldual}, for every $\bb\in\Int\mathscr M$ the minimum of the dual
equals the maximal entropy attained at the unique primal optimizer,
\[
S(\bb):=\min_{\by\in\R^m} f(\by)=-\Tr\bigl(X^\star\log X^\star\bigr),
\qquad X^\star=\exp_1 A(\by^\star).
\]
Figure~\ref{fig:entropy_landscape} displays this maximal von Neumann entropy over
$\mathscr M$. It is concave and peaks at $\log 3$ at the center $\bb=\mathbf 0$, where
$X^\star=\tfrac13 I$ is maximally mixed; it decreases to zero on the strictly convex
(elliptic) part of the boundary, where $X^\star$ is a rank-one projector (a pure state);
and it stays elevated, up to $\log 2$, along the two flat edges, where $X^\star$ has
rank two. The map $\bb\mapsto S(\bb)$ is exactly the negative Legendre conjugate of the
log-partition function $g(\by)=\log\Tr\exp A(\by)$, so its flattening near
$\partial\mathscr M$ is the dual face of the curvature estimates of
Section~\ref{curvature}: the problem stiffens precisely where $\bb$ approaches the boundary
and the smallest eigenvalue of $X^\star$ collapses.

\begin{figure}[h]
	\begin{center}
		\includegraphics[width=\textwidth]{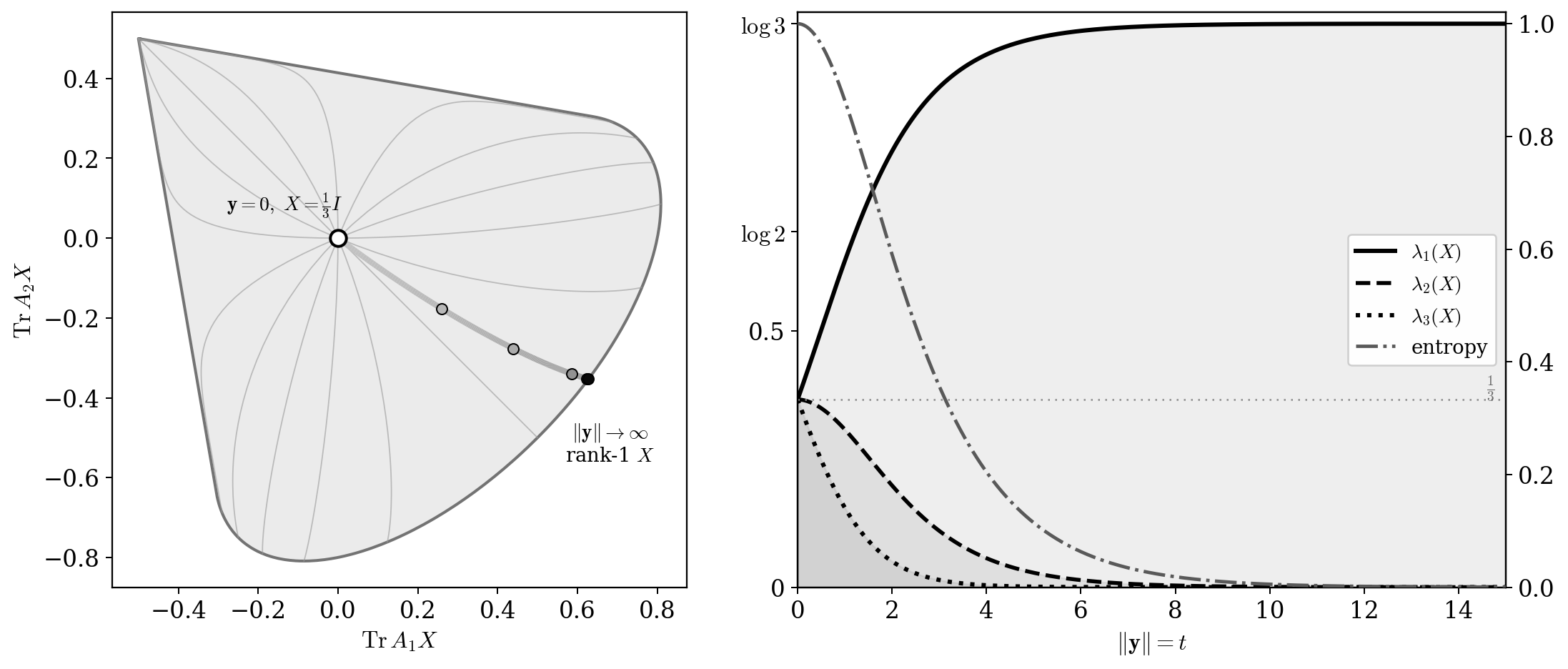}
		\caption{Gibbs state $X(\by)=\exp_1 A(\by)$ along the ray $\by=t\,\bv$.
			Left: its image $\mathcal A(X(\by))$ in $\mathscr M$ (gray curve darkening with
			$t$), from the center $X=\tfrac13 I$ ($t=0$, white marker) to the boundary as
			$t\to\infty$; faint gray curves are the images of other rays, filling
			$\Int\mathscr M$ as in Lemma~\ref{diffeo}. Right: the eigenvalues of $X(\by)$
			($\lambda_1$ solid, $\lambda_2$ dashed, $\lambda_3$ dotted, right axis) and the entropy $S$
			(dash-dotted, left axis). The spectrum spreads from
			$(\tfrac13,\tfrac13,\tfrac13)$ to $(1,0,0)$ while $S$ falls from $\log 3$ to
			$0$.\label{fig:gibbs_spectrum}}
	\end{center}
\end{figure}

Figure~\ref{fig:gibbs_spectrum} follows a single ray $\by=t\,\bv$, $t\ge 0$. On
the left, its image $\mathcal A(X(\by))$ under the forward map $\by\mapsto\nabla g(\by)$
runs from the center ($\by=0$, $X=\tfrac13 I$) to a boundary point as $t\to\infty$;
the pencil of rays through the origin maps to a fan of curves filling $\Int\mathscr M$,
in agreement with the diffeomorphism of Lemma~\ref{diffeo}. On the right, the
eigenvalues $\lambda_1\ge\lambda_2\ge\lambda_3$ of $X(\by)$ — partitioning the unity, since $\trace X(\by)=1$
 — spread from the flat profile $(\tfrac13,\tfrac13,\tfrac13)$ toward the
rank-one limit $(1,0,0)$, while the entropy decreases monotonically from $\log 3$ to
$0$. This degeneration is the mechanism behind the conditioning of $f$: as
$t\to\infty$ the iterate reaches $\partial\mathscr M$ and
$\lambda_{\min}(X(\by))\to 0$, so the strong-convexity constant of $f$, proportional to
$\lambda_{\min}(X)$ by Lemma~\ref{curvature}, vanishes — the geometric content of the
margin $\delta$ in the convergence analysis of Section~\ref{convergence}.

\subsection{Dense random instances}
{
	We first compared \texttt{maxentmom} with MOSEK and SDPNAL+ on dense (non-sparse) moment-body
	membership problems\footnote{The benchmark is available at
		\url{https://homepages.laas.fr/henrion/software/maxentmom/dense.m}.} generated as follows: $m$
	normally distributed random symmetric matrices $A_1,\dots,A_m$, a feasible point $X=\exp_1 Z$ with
	$Z$ a normally distributed random symmetric matrix, and $\bb=\mathcal A(X)\in\Int\mathscr M$. To
	expose the dependence on the two sizes separately, we vary $m$ and $n$ independently rather than only
	along the diagonal $m=n$. Our solver \texttt{maxentmom} is applied to the preconditioned data
	produced by Algorithm~\ref{algo}, with the preconditioning and the solve timed separately; MOSEK and
	SDPNAL+ are applied, without preconditioning, to the feasibility semidefinite program (zero objective
	function) with the constraints $\langle A_i,X\rangle=b_i$ and the explicit trace constraint
	$\trace X=1$. All solvers target accuracy $10^{-8}$. Timings (seconds) are reported in
	Table~\ref{tab:dense}.
	
	\begin{table}[ht]
		\centering
		\begin{tabular}{rr|cc|cc}
			& & \multicolumn{2}{c|}{\tt maxentmom} & & \\
			$m$ & $n$ & precond. & solve & MOSEK & SDPNAL+ \\
			\hline
			\multicolumn{6}{l}{\footnotesize fixed $n=100$, increasing $m$}\\
			$25$   & $100$ & $0.002$ & $0.005$ & $0.14$  & $0.012$\\
			$100$  & $100$ & $0.007$ & $0.013$ & $0.60$  & $0.033$\\
			$400$  & $100$ & $0.055$ & $0.071$ & $1.99$  & $0.122$\\
			$1500$ & $100$ & $0.904$ & $0.693$ & $27.18$ & $2.145$\\
			\hline
			\multicolumn{6}{l}{\footnotesize fixed $m=100$, increasing $n$}\\
			$100$ & $50$  & $0.003$ & $0.005$ & $0.14$  & $0.014$\\
			$100$ & $100$ & $0.007$ & $0.013$ & $0.60$  & $0.033$\\
			$100$ & $200$ & $0.015$ & $0.035$ & $2.04$  & $0.075$\\
			$100$ & $300$ & $0.062$ & $0.092$ & $6.01$  & $0.120$\\
			\hline
			\multicolumn{6}{l}{\footnotesize rectangular instances}\\
			$200$ & $100$ & $0.014$ & $0.017$ & $1.06$  & $0.033$\\
			$300$ & $150$ & $0.116$ & $0.043$ & $3.70$  & $0.113$\\
			$400$ & $200$ & $0.204$ & $0.087$ & $10.53$ & $0.243$\\
			$150$ & $300$ & $0.128$ & $0.080$ & $8.48$  & $0.206$\\
			$200$ & $400$ & $0.331$ & $0.148$ & $22.29$ & $0.400$\\
			$500$ & $250$ & $0.523$ & $0.222$ & $25.22$ & $0.485$\\
			\hline
		\end{tabular}
		\caption{Computing time in seconds on dense random moment-body membership problems, for varying
			constraint number $m$ and matrix size $n$. For \texttt{maxentmom} the preconditioning
			(Algorithm~\ref{algo}) and the L-BFGS solve are timed separately; MOSEK and SDPNAL+ receive the
			original, non-preconditioned data. MOSEK and SDPNAL+ reach a normalized residual below $10^{-11}$ on
			every instance; \texttt{maxentmom} reaches a few $10^{-9}$ on all instances except $m=1500$, $n=100$
			(see text).}
		\label{tab:dense}
	\end{table}
	
	Three points stand out. First, on all but one instance preconditioning does what the curvature
	analysis of Sections~5--7 predicts: L-BFGS reaches a normalized residual of a few $10^{-9}$ in of the
	order of ten iterations, so the solve is little more than a handful of eigenvalue decompositions. The
	exception is the constraint-heavy instance $m=1500$, $n=100$ (about $30\%$ of the traceless
	dimension), where the residual stalls at $4\cdot10^{-4}$: there the minimizer has
	$\lambda_{\min}(X^*)\approx2\cdot10^{-14}$, so $\bb$ lies essentially on $\partial\mathscr M$, and the
	margin $\delta$ is minuscule. This is precisely the ill-conditioned regime predicted by the bound
	$\kappa\le\tfrac12\,n^{1+\sqrt2/\delta}$, compounded by the near-singularity of the
	centered Gram matrix that the whitening must invert; MOSEK and SDPNAL+, working on the original data,
	are unaffected and still reach $10^{-12}$.
	
	Second, varying $m$ and $n$ separately exposes the two-part cost of \texttt{maxentmom}. The solve
	grows with $n$, through the $O(n^3)$ eigendecomposition at each iteration (and the comparable
	$O(mn^2)$ assembly of $A(\by)$ and of the gradient at moderate $n$), and only mildly with $m$; the
	preconditioning, whose cost is the $O(m^2n^2)$ formation of the centered Gram matrix together with
	the $O(m^3)$ whitening, grows roughly quadratically in each of $m$ and $n$, the cubic term taking over
	at the largest $m$ (the jump to $0.9$\,s at $m=1500$). Because the solve is so cheap, the
	preconditioning becomes the larger part of \texttt{maxentmom}'s total at the larger sizes --- about
	$70\%$ at $n\ge200$ or $m\ge400$.
	
	Third, the comparison with the two reference solvers is consistent across aspect ratios. The
	\texttt{maxentmom} solve alone is two to three times faster than SDPNAL+ on every instance; once
	preconditioning is included, the total is faster than SDPNAL+ on the smaller and tall ($m>n$)
	instances and within a factor of about $1.5$ on the largest and wide ($n>m$) ones, where the
	$O(m^2n^2+m^3)$ preconditioning dominates. Both first-order methods vastly outperform the
	second-order interior-point solver MOSEK, by one to two orders of magnitude: at $n=400$
	($m=200$), \texttt{maxentmom} and SDPNAL+ are about $47$ and $56$ times faster respectively, and the
	gap widens with $n$, consistently with the steep growth of the interior-point Schur-complement
	computations.
}

{
\subsection{Matrix completion problems}

A natural family of membership oracles arises from positive semidefinite matrix
completion. Let $\Omega\subseteq\{(k,l):1\le k\le l\le n\}$ be a set of matrix positions
containing the whole diagonal, and to each position associate the orthonormal symmetric
matrix
\[
A_{kk}=e_ke_k^T,\qquad A_{kl}=\tfrac{1}{\sqrt2}(e_ke_l^T+e_le_k^T)\quad k<l,
\]
so that $\langle A_i,A_j\rangle=\delta_{ij}$ and the map $\mathcal A$ simply reads off the
sampled entries of a matrix.
Deciding whether a matrix specified on $\Omega$ admits a
positive semidefinite, unit trace completion is then exactly the membership problem
$\bb\in\mathscr M=\mathcal A(\S^n_1)$. Note that in this original matrix-completion formulation, Assumption~\ref{injective} is not literally satisfied. Indeed, since the whole diagonal is sampled, the diagonal measurement matrices satisfy $\sum_{k=1}^n A_{kk}=I_n$,
so the identity belongs to the span of the measurements. This only reflects a harmless redundancy: the trace-one constraint is already encoded by the sampled diagonal entries, and the dual log-partition function is flat in the corresponding identity direction, while the primal matrix $X=\exp_1 A(\by)$ is unchanged by adding a scalar multiple of $I_n$ to $A(\by)$. Equivalently, one may quotient out this redundant direction, or replace the $n$ diagonal measurements by $n-1$ traceless diagonal contrasts. The resulting reduced formulation is equivalent to the original completion problem and satisfies the full-dimensionality Assumption~\ref{injective}. In contrast to the
dense projections of the previous experiments, the matrices $A_i$ are now extremely sparse,
so the data matrix $A$ has only $O(|\Omega|)$ nonzero entries and each gradient evaluation
$\mathcal A(X(\by))$ is cheap; matrix completion thus exercises the complementary,
sparse regime of the oracle.

We generate random feasible instances as follows.\footnote{The benchmark is available at  \url{https://homepages.laas.fr/henrion/software/maxentmom/mc.m}.}
We draw a Gaussian matrix $W\in\R^{n\times n}$ and set $X_0=WW^T/\trace(WW^T)$, a full rank
density matrix, so that $\bb=\mathcal A(X_0)\in\Int\mathscr M$. The set $\Omega$ consists of
the full diagonal together with a fraction $p$ of the off-diagonal positions drawn uniformly
at random, with $p\in\{2,5,10,20\}\%$ and $n\in\{100,200,500,1000,2000,5000\}$; for the
largest instances the number of affine constraints $m=|\Omega|$ reaches about
$2.5\times10^6$. On each instance we run our  solver \texttt{maxentmom} and, for
comparison, SDPNAL+ applied to the feasibility semidefinite program (with zero objective function) with the same constraints, both to a residual tolerance $10^{-8}$. The computing
times are collected in Table~\ref{tab:completion}.

\begin{table}[ht]
	\centering
	\begin{tabular}{r|cccc}
		\hline
		& \multicolumn{4}{c}{sampling density $p$}\\
		$n$ & $2\%$ & $5\%$ & $10\%$ & $20\%$\\
		\hline
		$100$  & $0.02/0.30$    & $0.01/0.05$    & $0.01/0.04$    & $0.01/0.02$\\
		$200$  & $0.02/0.16$    & $0.03/0.04$    & $0.03/0.02$    & $0.03/0.02$\\
		$500$  & $0.13/0.10$    & $0.16/0.08$    & $0.15/0.09$    & $0.21/0.08$\\
		$1000$ & $0.66/0.40$    & $0.72/0.37$    & $0.76/0.39$    & $1.04/0.47$\\
		$2000$ & $4.98/6.01$    & $9.09/4.85$    & $7.97/4.17$    & $11.54/4.92$\\
		$5000$ & $112.85/32.69$ & $151.93/50.37$ & $168.01/29.84$ & $171.53/31.54$\\
		\hline
	\end{tabular}
	\caption{Computing time in seconds for the matrix completion membership oracle, reported as
		\texttt{maxentmom}\,/\,SDPNAL+, as a function of the matrix size $n$ and the sampling density
		$p$ (the fraction of off-diagonal entries revealed, in addition to the full diagonal).}
	\label{tab:completion}
\end{table}

The two solvers stay within a small factor of one another over the whole range, which is
notable since \texttt{maxentmom} is a rudimentary prototype whereas SDPNAL+ is a mature,
compiled code. For small problems the absolute times are dominated by fixed overheads --- the
$0.30$ s of SDPNAL+ at $n=100$, $p=2\%$ is its start-up cost --- and \texttt{maxentmom} is
competitive or faster: it is faster at every density for $n=100$, and faster at low density
for $n=200$. In the intermediate range $500\le n\le 2000$ the two remain within a factor of
about $2.5$, with \texttt{maxentmom} still the faster of the two at the sparsest sampling,
where it wins even at $n=2000$ ($4.98$ against $6.01$ s at $2\%$). At $n=5000$ the balance has
tipped to SDPNAL+, which is then faster by a factor of three to five.

The two methods respond very differently to the sampling density. The cost of
\texttt{maxentmom} increases with $p$ --- at $n=5000$ it rises from $113$ to $172$ seconds
between $2\%$ and $20\%$ --- since a larger constraint set raises both the dimension of the
dual variable and the number of L-BFGS iterations. SDPNAL+, by contrast, is essentially
insensitive to $p$ (between $30$ and $50$ seconds at $n=5000$, irrespective of density): its
running time is governed by the eigenvalue decompositions in the projection onto the
semidefinite cone, which scale with $n$ but not with $m$, the affine projection being trivial
here because $\mathcal A\mathcal A^T$ is diagonal. Both methods scale roughly cubically in
$n$, consistent with the cubic complexity established above and with the cost being dominated
by dense eigenvalue computations; \texttt{maxentmom} shows a mild super-cubic growth at fixed
density, reflecting the slow increase of its iteration count with $n$.
}

\section{Conclusion} 

Motivated by pre-conditioning strategies for semidefinite optimization, this paper reports on a specific problem class whose geometry is simple enough to allow for a comprehensive analysis. We consider the moment body membership oracle problem, which consists of determining whether a given vector of size $m$ belongs to a given linear projection of the spectraplex, the compact convex set of unit trace positive semidefinite matrices of size $n$-by-$n$. Inspired by maximum entropy techniques from quantum information theory, we propose to solve the problem by minimizing on the whole $m$-dimensional space a dual smooth strictly convex log-partition function. Geometric curvature analysis reveals how key input data quantities can be modified to improve the problem conditioning. After pre-conditioning, we can solve the convex dual problem with L-BFGS, a widely used first-order algorithm approximating second-order information with limited gradient evaluation and storage. 

Numerical experiments with a rudimentary Matlab implementation, \texttt{maxentmom}, support the analysis in two complementary regimes. On dense random instances, the L-BFGS solve is very fast once the data have been preconditioned, often requiring only a few iterations; the main cost then becomes the preliminary centering and whitening step. The method is substantially faster than MOSEK on these tests and remains competitive with SDPNAL+, although the latter can be faster once the preconditioning cost is included. The experiments also confirm the role of the geometric margin: well-conditioned interior points are solved accurately and quickly, whereas instances close to the boundary of the moment body are much harder. On sparse matrix-completion problems, where the measurements are simple coordinate observations, \texttt{maxentmom} and SDPNAL+ have comparable behavior over a wide range of sizes, with both methods ultimately limited by dense spectral computations.

A natural direction for future work is to reduce this spectral bottleneck. Instead of computing the full matrix exponential and its eigendecomposition, one could estimate the traces defining the objective and gradient by randomized probing, using only matrix-vector products with the exponential. This idea is closely related to recent randomized entropic methods \cite{l23,cl25b}. It is especially promising for genuinely sparse problems, although the stochastic estimates introduce a bias-variance tradeoff and may be better suited to moderate accuracy or as a preliminary phase before a final deterministic refinement. Understanding when such randomized variants outperform exact spectral computations remains an important open question.

Polynomial SOS decompositions are particular cases of the moment body membership oracle, where $X$ is the Gram matrix representing a polynomial $p(\bx) = \phi^T(\bx) X \phi(\bx)$ as a quadratic form w.r.t. some basis vector $\phi$. The linear map $\mathcal A(X)=\bb$ matches $X$ with $p$ expressed as a coefficient vector $\bb$ in some basis.
It can be normalized since $\int p(\bx) d\mu(\bx) = \trace (X\int \phi(\bx) \phi^T(\bx) d\mu(\bx)) = \trace X$
whenever $\phi$ is an orthonormal basis with respect to the inner product induced by $\mu$. Memberships in truncated quadratic modules, also called weighted SOS decompositions, can also be modeled as particular moment body problems. They are at the core of the moment-SOS hierarchy for polynomial optimization \cite{hkl20}. It would be interesting to derive specific curvature properties to pre-condition these problems in the same way we did it for general moment bodies. Relationships with the dual certificates of truncated quadratic module membership investigated in \cite{dp22} are also worth investigating, especially since these dual certificates allow to construct SOS representations with rational coefficients.

In the context of semidefinite relaxations of combinatorial optimization problems, the trace one constraint holds for the first relaxation of the moment-SOS hierarchy. This constant trace property was exploited in \cite{hr00} in the context of spectral bundle methods. It was generalized in \cite{mlm23,mlmw22} where it was shown that every polynomial optimization problem on a compact semialgebraic set has an equivalent equality constrained formulation on a sphere (possibly after adding some artificial variables), and hence a constant trace moment relaxation.

A natural extension of our approach consists of minimizing a linear function on the moment body, i.e. given a matrix $C \in \S^n$, solving the semidefinite optimization problem
$$
\min_{X \in \S^n_1} \trace(CX) \quad\mathrm{s.t.}\quad \mathcal A(X) = \bb. 
$$
For a given regularization parameter $\mu>0$, to a primal entropic problem
$$
\min_{X \in \S^n_1} \trace(CX) - \mu\:\trace(X - X \log X)  \quad\mathrm{s.t.}\quad \mathcal A(X) = \bb
$$
corresponds a dual log-partition problem
$$ 
\max_{\by\in\R^m} \bb^T \by - \mu\;\log \trace \exp  (-\tfrac{1}{\mu}(C - A(\by))).
$$
One then follows a primal admissible central path
$$
X^*_{\mu} = \exp_1(-\frac{1}{\mu}(C-A(\by^*_{\mu}))) \in \mathcal A^{-1}(\bb)
$$
parametrized by dual optimal solutions $\by^*_{\mu}$
and we let $\mu \to 0^+$. A detailed analysis of convergence of this method remains to be done. Note that the idea was followed recently in \cite{l23,cl25}, but without the trace-one restriction. Consequently, the dual function there is the much less regular partition function, which is the exponential of the log-partition function. This may explain why the semidefinite optimization experiments reported in \cite{cl25} are somewhat disappointing. Whether more convincing and scalable numerical results can be obtained with the log-partition function remains however to be seen.

\section*{Acknowledgement}

Solving the spectrahedral shadow membership with a first order optimization algorithm was suggested to the author by Stephan Weis at Mathematisches Forschungsinstitut Oberwolfach in August 2024. He also pointed out references \cite{h24,jv23} {and noticed several mistakes in the first version of this paper.} A significant part of this work was done during a stay at the Institute of Pure and Applied Mathematics of the University of California at Los Angeles, whose hospitality has been appreciated.
This work benefited from feedback from Saroj Prasad Chhatoi, Jean Bernard Lasserre, Victor Magron, Jiawang Nie {as well as Samuel Burer and two anonymous reviewers}. {The author acknowledges the use of AI for assistance with brainstorming idea, mathematical development, coding and drafting the manuscript. The final content, analysis and conclusions remain the sole responsibility of the author.}

\end{document}